\documentclass[12pt,reqno]{amsart}
\usepackage{amssymb,delarray}
\usepackage{graphicx}

\usepackage{epsf,amsfonts,amscd,latexsym}

\def\leq{\leqslant}
\def\geq{\geqslant}

\def\edvo{\rule {6pt}{6pt}} 

\newtheorem{thm}{Theorem}[section]
\newtheorem{lem}[thm]{Lemma}
\newtheorem{prop}[thm]{Proposition}

\newtheorem{cor}[thm]{Corollary}

\newtheorem{que}[thm]{Question}

{\catcode`\@=11
\gdef\n@te#1#2{\leavevmode\vadjust{%
 {\setbox\z@\hbox to\z@{\strut#1}%
  \setbox\z@\hbox{\raise\dp\strutbox\box\z@}\ht\z@=\z@\dp\z@=\z@%
  #2\box\z@}}}
\gdef\leftnote#1{\n@te{\hss#1\quad}{}}
\gdef\rightnote#1{\n@te{\quad\kern-\leftskip#1\hss}{\moveright\hsize}}
\gdef\?{\FN@\qumark}
\gdef\qumark{\ifx\next"\DN@"##1"{\leftnote{\rm##1}}\else
 \DN@{\leftnote{\rm??}}\fi{\rm??}\next@}}

\begin{document}
\baselineskip=14.pt plus 2pt 

\title[
] {On 
Symplectic  Coverings of the Projective Plane}
\author[G.-M.~Greuel]{G.-M.~Greuel}
\address{Kaiserslautern University}
 \email{greuel@mathematik.uni-kl.de}

 \author[Vik.S.~Kulikov]{Vik.S.~Kulikov}
\address{Steklov Mathematical Institute}
 \email{kulikov@mi.ras.ru}
\dedicatory{} \subjclass{}
\thanks{The first author was partially supported by the  DFG-Schwerpunkt
  ''Globale Methoden in der komplexen Geometrie'' and the second by the RFBR 
(
02-01-00786) and by the DFG (436 RUS 17/84/03).} 
\begin{abstract}
We prove that a resolution of singularities of any finite covering
of the projective plane branched along a Hurwitz curve $\bar H$
and, maybe, along a line "at infinity"  can be embedded as a
symplectic submanifold into some projective algebraic manifold
equipped with an integer K\"{a}hler symplectic form (assuming that
if $\bar H$ has negative nodes, then the covering is non-singular
over them). For cyclic coverings we can realize this embeddings into a
rational algebraic $3$--fold. Properties of the Alexander polynomial of $\bar{H}$ are investigated and applied to the calculation of
the first Betti number $b_1(\overline X_n)$ of a resolution
$\overline X_n$ of singularities of  $n$-sheeted cyclic coverings
of $\mathbb C\mathbb P^2$ branched along $\bar H$
and, maybe, along a line "at infinity".  We prove that
$b_1(\overline X_n)$ is even if $\bar H$ is an irreducible Hurwitz
curve but, in contrast to the algebraic case, that it can take any non-negative value in the case when
$\bar H$ consists of several irreducible components.
\end{abstract}

\maketitle
\setcounter{tocdepth}{2}


\def\st{{\sf st}}

\setcounter{section}{-1}
\section{Introduction}

The notion of Hurwitz curves in the projective plane $\mathbb
C\mathbb P^2$ with respect to a linear projection $\text{pr}:\mathbb C\mathbb P^2\smallsetminus\{p_\infty\}\to\mathbb C\mathbb
P^1, (p_\infty$ denotes the centre of the projection $\text{pr}$)
was introduced in \cite{Moi2} and is a natural generalization of the
notion of plane algebraic curves (in \cite{Moi2}, Hurwitz curves
are called "semi-algebraic" curves). A precise definition of
Hurwitz curves can be found, for example, in \cite{Kh-Ku}. In this paper
we give another equivalent (see Lemma \ref{definition}) definition
of Hurwitz curves as follows. Let $\mathbb C^2_i$ be two copies of
the affine plane $\mathbb C^2$, $i=1,2$, with coordinates
$(u_i,v_i)$, $u_2=1/u_1$ and $v_2=v_1/u_1$, which cover $\mathbb
C\mathbb P^2\setminus p_\infty$,
 such that $\text{pr}$ is given by $(u_i,v_i)\to u_i$
in the charts $\mathbb C^2_i$. A set $\bar H\subset \mathbb
C\mathbb P^2\setminus \{p_\infty\}$, closed in $\mathbb C\mathbb P^2$, is called a {\it
Hurwitz curve of degree} $m$ if, for $i=1,2$, $\bar H\cap \mathbb C^2_i$
coincides with the set of zeros of an equation
\begin{equation} \label{function} F_i(u_i,v_i):=v_i^m
+\sum_{j=0}^{m-1}c_{j,i}(u_i)v_i^j=0
\end{equation}
such that
\begin{itemize}
\item[($i$)] $F_i (u_i, v_i)$ is a  $C^\infty$--smooth complex valued function in
$\mathbb C^2$; \item[($ii$)] the function $F_i(u_i,v_i)$ has only
a finite number of critical values, that is, there are finitely
many values of $u_i$, say $u_{i,1},\dots, u_{i,n_i}$, such that
the polynomial equation
\begin{equation} \label{def} v_i^m
+\sum_{j=0}^{m-1}c_{j,i}(u_{i,0})v_i^j=0\end{equation} has no
multiple roots for $u_{i,0}\not\in \{ u_{i,1},\dots, u_{i,n_i}
\}$; \item[($iii$)] if $v_{i,j}$ is a multiple root of equation
(\ref{def}) for $u_{i,j}\in \{ u_{i,1},\dots, u_{i,n_i} \}$, then,
in a neighbourhood of the point $(u_{i,j},v_{i,j})$ (which we call
a {\it critical point} of $\bar H$), the set $\bar H$ coincides
with the solution a complex analytic equation.
\end{itemize}

Note that after rescaling $\widetilde v_i=\varepsilon v_i$,
$0<\varepsilon <<1$, Hurwitz curves become symplectic surfaces in
$\mathbb C\mathbb P^2$ (see also proof of Theorem 3.1).

More general, one can consider so called {\it topological Hurwitz
curves} which have {\it cone singularities} (see the definition of
cone singularities in \cite{Kh-Ku}).

A Hurwitz (resp. topological Hurwitz) curve $\bar H$ is called
{\it irreducible} if $\bar H\setminus M$ is connected for any
finite set $M\subset \bar H$, and we say that a Hurwitz curve
$\bar H$ {\it consists of $k$ irreducible components} if $$k=\max
\# \{ \text{connected components of}\, \bar H\setminus M\},$$
where the maximum is taken over all finite sets $M\subset \bar H$.

Let $H$ be an {\it affine Hurwitz curve}, that is, $H=\bar H\cap
(\mathbb C\mathbb P^2\setminus L_{\infty})$, where
$L_{\infty}$ is a line which is a fibre of $\text{pr}$ 
being in general position with respect to $\bar H$. Then the
fundamental group $\pi_1=\pi_1(\mathbb C\mathbb P^2\setminus (\bar
H\cup L_{\infty}))$ does not depend on the choice of $L_{\infty}$
and belongs to the class $\mathcal C$ of so called $C$-groups.

By definition, a {\it $C$-group} is a group together with  a
finite presentation
\begin{equation} \label{zero}
G_W=<x_1,\dots ,x_m \, \mid \, x_i= w_{i,j,k}^{-1} x_jw_{i,j,k} ,
\, \, w_{i,j,k}\in W\, >,
\end{equation}
where $W=\{ w_{i,j,k}\in \mathbb F_m\, \mid \, 1\leq i,j\leq m,\,
\, 1\leq k\leq h(i,j)\}$ is a subset of elements of the free group
$\mathbb F_m$ (it is possible that
$w_{i_1,j_1,k_1}=w_{i_2,j_2,k_2}$ for $(i_1,j_1,k_1)\neq
(i_2,j_2,k_2)$), generated by free generators $x_1,\dots,x_m$ and
$h:\{1,\dots , m\}^2 \to \mathbb Z$ is some function. Such a
presentation is called a $C$-{\it presentation} ($C$, since all
relations are conjugations). Let $\varphi_W: \mathbb F_m\to G_W$
be the canonical epimorphism. The elements $\varphi_{W}(x_i)\in
G$, $1\leq i\leq m$, and the elements conjugated to them are
called the {\it $C$-generators} of the $C$-group $G$. Let
$f:G_1\to G_2$ be a homomorphism of $C$-groups. It is called a
{\it $C$-homomorphism} if the images of the $C$-generators of
$G_1$ under $f$ are $C$-generators of the $C$-group $G_2$.
$C$-groups will be considered up to $C$-isomorphisms.

A $C$-presentation (\ref{zero}) is called a {\it Hurwitz
$C$-presentation of degree} $m$ if for each $i=1,\dots,m$ the word
$w_{i,i,1}$ coincides with the product $x_1\dots x_m$, and a
$C$-group $G$ is called a {\it Hurwitz $C$-group of degree} $m$ if
it possesses a Hurwitz $C$-presentation of degree $m$. In other
words, a $C$-group $G$ is a Hurwitz $C$-group of degree $m$ if
there are $C$-generators $x_1,\dots, x_m$ generating $G$ such that
the product $x_1\dots x_m$ belongs to the center of $G$. Note that
the degree of a Hurwitz $C$-group $G$ is not defined canonically
and depends on the Hurwitz $C$-presentation of $G$. Denote by
$\mathcal H$ the class of all  Hurwitz $C$-groups.

Let $\bar H$ be a Hurwitz (resp. topological Hurwitz) curve of
degree $m$. A Zariski -- van Kampen presentation of
$\pi_1=\pi_1(\mathbb C^2\setminus H)$ (where $\mathbb C^2=\mathbb
P^2\setminus L_{\infty}$ and $L_{\infty}$, a fibre of pr, is in general
position with respect to $\bar H$) defines on $\pi_1$ a structure
of a Hurwitz $C$-group of degree $m$ (see \cite{Ku}), and in
\cite{Ku}, it was proved that any Hurwitz $C$-group $G$ of degree $m$
can be realized as the fundamental group $\pi_1(\mathbb
C^2\setminus H)$ for some Hurwitz curve $\bar H$ with
singularities of the form $w^m-z^m=0$, $\deg \bar H=2^nm$, where
$n$ depends on the Hurwitz $C$-presentation of $G$. Since we
consider $C$-groups up to $C$-isomorphisms, the class
$\mathcal H$ coincides with the class $\{ \, \pi_1(\mathbb
C^2\setminus H)\, \}$ of fundamental groups of the complements of
affine Hurwitz (resp. topological Hurwitz) curves.\\

A free group $\mathbb F_n$ with fixed free generators is a
$C$-group and for any $C$-group $G$ the canonical $C$-epimorphism
$\nu : G\to \mathbb F_1$, sending the $C$-generators of $G$ to the
$C$-generator of $\mathbb F_1$, is well defined. Denote by $N$
its kernel. Note that if all $C$-generators of a $C$-group $G$ are
conjugated to each other (such $C$-group is called {\it
irreducible}), then $N$ coincides with
\footnote{For a group $G$ we use the standard notation $G'$ for
its commutator subgroup and $G^{\prime\prime}$ for the
commutator subgroup of $G'$.} $G'$.

Let $G$ be a $C$-group. The $C$-epimorphism $\nu$ induces the
following exact sequence of groups
$$1\to N/N'\to G/N'\buildrel{\nu_*}\over\longrightarrow  \mathbb F_1\to 1.$$
The $C$-generator of $\mathbb F_1$ acts on $N/N'$ by conjugation
$\widetilde x^{-1} b\widetilde{x}$, where $n\in N$ and
$\widetilde{x}$ is one
  of the $C$--generators of $G$. Denote by $h$
this action and by $h_{\mathbb C}$ the induced action on
$N/N'\otimes \mathbb C$. The characteristic polynomial $\Delta
(t)=\det (h_{\mathbb C}-t\text{Id})$ is called the {\it Alexander
polynomial} of the $C$-group $G$ (if the vector space $N/N'\otimes
\mathbb C$ over $\mathbb C$ is infinite dimensional, then, by
definition, the Alexander polynomial $\Delta (t)\equiv 0$). For a
(topological) Hurwitz curve $\bar H$ the Alexander polynomial
$\Delta (t)$ of the group $\pi_1=\pi_1(\mathbb C^2\setminus H)$ is
called the {\it Alexander polynomial} of $\bar H$. Note that the
Alexander polynomial $\Delta (t)$ of a (topological) Hurwitz curve
$\bar H$ does not depend on the choice of the generic line $L_{\infty}$.

Let $G=<x_1,\dots,x_m\, \mid\, r_1,\dots,r_n>$ be a
$C$-presentation of a $C$-group $G$ and $\mathbb F_m$ be the free
group freely generated by the $C$-generators $x_1,\dots,x_m$.
Denote by $\frac{\partial}{\partial x_i}$ the Fox derivative
(\cite{C-F}), that is, the endomorphism of the group ring $\mathbb
Z[\mathbb F_m]$ over $\mathbb Z$ of the free group $\mathbb F_m$
into itself, such that $\frac{\partial}{\partial x_i} :\mathbb
Z[\mathbb F_m]\to \mathbb Z[\mathbb F_m]$ is a $\mathbb Z$-linear
map defined by the following properties
\begin{equation} \label{fox}
 \begin{array}{rcl}
\displaystyle \frac{\partial x_j}{\partial x_i} & = & \delta_{i,j} \\
\displaystyle \frac{\partial uv}{\partial x_i} & = & \displaystyle
\frac{\partial u}{\partial x_i}+u\frac{\partial v}{\partial x_i}
\end{array}
\end{equation}
for any $u,v\in \mathbb Z[\mathbb F_m]$. It is well known (it is
proved, for example, in \cite{Ro} in the case of knot groups and
generalized to the case of $C$-groups in \cite{Ku3}) that the
greatest common divisor of the minors of order $m-1$ in the
matrix 
$$\nu_*\bigl(\tfrac{\partial r_i}{\partial x_j}\bigr)\in \text{Mat}_{n\times
m}(\mathbb Z[t,t^{-1}])$$  coincides with the Alexander polynomial
$\Delta (t)$ of $G$ up to a factor $\pm t^k$ invertible in
$\mathbb Z[t,t^{-1}]$, where $r_i$, $i=1,\dots,n$, are the
defining relations of $G$ and $\nu_*: \mathbb Z[\mathbb F_m]\to
\mathbb Z[\mathbb F_1]\simeq \mathbb Z[t,t^{-1}]$ is induced by
the canonical $C$-epimorphism $\nu :\mathbb F_m\to \mathbb F_1$.

The properties of the Alexander polynomials of  plane algebraic
curves and their application to the calculation of the first Betti
number of a cyclic covering of the projective plane are well-known
(see, for example, \cite{Zar}, \cite{Lib}, \cite{R}, \cite{E},
\cite{Ko}, \cite{Ku3}, \cite{Ku-Ku}). One of the aims of this article is
to generalize these results to the case of Hurwitz curves and to
apply them to the calculation of the first Betti number of a
cyclic covering of the projective plane branched along a Hurwitz
curve.

The main results of this article are the following theorems and
corollaries.
\begin{thm} \label{main} Let $\bar H$ be a (topolgical) Hurwitz curve of degree $d$
and $\Delta (t)$ its Alexander polynomial. Then
\begin{itemize}
\item[($i$)] $\Delta (t)\in \mathbb Z[t]$; \item[($ii$)]
$\Delta(0)=\pm 1$; \item[($iii$)] the roots of $\Delta (t)$  are
$d$-th roots of unity; \item[($iv$)] the action of $h_{\mathbb C}$
on $(N/N')\otimes \mathbb C$ is semisimple.
\end{itemize}
\end{thm}
Moreover, the Alexander polynomial $\Delta (t)$ of a Hurwitz curve
$\bar H$ of degree $d$ is a divisor of the polynomial
$(t-1)(t^d-1)^{d-2}$ (see Theorem \ref{div}) and if $\bar H$
consists of $k$ irreducible components, then the multiplicity of
the root $t=1$ of its Alexander polynomial $\Delta (t)$ is equal
to $k-1$ (see Theorem \ref{mult}).
\begin{thm}\label{Roots}
If $\bar H$ is an irreducible (topological) Hurwitz curve, then
\begin{itemize}
\item[($i$)] $\Delta (t)$ is a reciprocal polynomial, i.e.
$\Delta(t)=t^{\deg \Delta(t)} \Delta(t^{-1})$; \item[($ii$)] $\deg
\Delta (t)$ is an even number; \item[($iii$)] $\Delta (1)=1$.
\end{itemize}
\end{thm}

\begin{cor} \label{cor} Let $\bar H$ be an irreducible (topological) Hurwitz curve
of $\deg \bar H=p^n$, where $p$ is a prime number. Then
\begin{itemize}
\item[($i$)] $\Delta (t)\equiv 1$; \item[($ii$)]the group
$\pi_1^{\prime}/\pi_1^{\prime \prime}$ is a finite group, where
$\pi_1=\pi_1(\mathbb C^2\setminus H)$.
\end{itemize}
\end{cor}

Note also that if $J$ is an almost complex structure in $\mathbb
C\mathbb P^2$ compatible with the Fubini -- Studi symplectic form
and if $\bar H$ is a $J$-holomorphic curve in $\mathbb C\mathbb P^2$
of degree $m$, that is, the class $[\bar H] \text{ equals }  m[\mathbb C\mathbb
P^1$] in $H_2(\mathbb C\mathbb P^2,\mathbb Z)$, then
$\pi_1=\pi_1(\mathbb C\mathbb P^2\setminus (\bar H\cup
L_{\infty}))$ is a Hurwitz $C$-group of degree $m$, where
$L_{\infty}$ is one of the $J$-lines being in general position
with respect to $\bar H$. Indeed, if we chose a pencil of
pseudo-holomorphic lines having $L_{\infty}$ as a member, then, by
the Zariski--van Kampen Theorem, a presentation of $\pi_1$ is
defined by a braid monodromy factorization of $\bar H$ with
respect to the chosen pencil. Therefore $\pi_1$ is a $C$-group,
and similar to the case of Hurwitz curves, it is easy to show (see
the proof of Theorem 6.1 from \cite{Ku}) that it is a Hurwitz
$C$-group of degree $m$. Thus, the Alexander polynomial of
a pseudo-holomorphic curve can be defined similarly and has
the same properties as in the case of Hurwitz curves.\\

The homomorphism $\nu :\pi_1\to \mathbb F_1$, where
$\pi_1=\pi_1(\mathbb C^2\setminus H)$ is the fundamental group of
the complement of an affine Hurwitz curve, defines an infinite
unramified cyclic covering $f=f_{\infty}:X'_{\infty}\to X'=\mathbb
C^2\setminus H$. We have $H_1(X'_{\infty},\mathbb Z)=N/N^{\prime}$
and the action of $h$ on $H_1(X'_{\infty},\mathbb Z)$ coincides
with the action of a generator of the covering transformation
group of the covering $f_{\infty}$. As it follows from \cite{Ku.O},
the\ group\ $H_1(X'_{\infty},\mathbb Z)$ is finitely generated. For
any $n\in \mathbb N$ denote\ by\! $\text{mod}_{n}:\mathbb F_1\to
\mu_n=\mathbb F_1/\{h^n\}$ the\ na\-tu\-ral epimorphism to the cyclic
group $\mu_n$ of degree $n$. The covering $f_{\infty}$ can be
factorized through the cyclic covering $f_n:X'_n\to \mathbb
C^2\setminus H$ associated with the epimorphism\ $\text{mod}_n\circ \nu$,
$f_{\infty}=g_n\circ f_n$. Since a Hur\-witz curve $\bar H$ has only
analytic singularities, the covering $f_n$ can be extended to a
smooth map ${\overline f}_n:{\overline X}_n\to \mathbb C\mathbb
P^2$ branched along $\bar H$ and, maybe, along $L_{\infty}$ (if
$n$ is not a divisor of $\deg \bar H$, then ${\overline f}_n$ is
branched along $L_{\infty}$), where $\overline X_n$ is a smooth
4-fold. The action $h$ induces an action $\overline h_n$ on
$\overline X_{n}$ and an action $\overline h_{n*}$ on
$H_1(\overline X_{n},\mathbb Z)$.

In section \ref{section4},  we show (see Theorem \ref{sympl-thm})
that any such $\overline X_n$ can be embedded as a symplectic
submanifold to a projective rational 3-fold on which the symplectic structure  is given by an integer K\"{a}hler form.

\begin{thm} \label{th} Let $\overline X_n$ be a resolution of singularities  of an $n$-sheeted cyclic covering
branched along a Hurwitz curve $\bar H$ and, maybe, along
$L_{\infty}$ and associated with the epimorphism $\text{mod}_n\circ \nu
:\pi_1\to \mathbb Z/n\mathbb Z$. Then the first Betti number
$$b_1(\overline X_n)=\dim_{\mathbb C} H_1(\overline X_n,\mathbb
C)=r_{n,\neq 1},$$ where $r_{n,\neq 1}$ is the number of roots of
the Alexander polynomial $\Delta (t)$ of the curve $\bar H$ which
are $n$-th roots of unity not equal to 1.
\end{thm}

Theorems \ref{main}, \ref{Roots}, \ref{th} and Corollary \ref{cor}
imply the following corollaries.
\begin{cor} Let $\overline X_n$ be a resolution of singularities of an $n$-sheeted cyclic covering
branched along a Hurwitz curve $\bar H$ and, maybe, along
$L_{\infty}$. If $\deg \bar H$ and $n$ are coprime, then
$b_1(\overline X_n)=0$.
\end{cor}
\begin{cor} Let $\overline X_n$ be a resolution of singularities of an $n$-sheeted cyclic covering
branched along an irreducible Hurwitz curve $\bar H$ and, maybe,
along $L_{\infty}$. Then $b_1(\overline X_n)$ is an even number.
\end{cor}
Moreover, we show that for any $k\in \mathbb N$, there is an
irreducible Hurwitz curve $\bar H_{k}$ such that for some $n$ (for
example, one can take $n=6$, see Proposition \ref{b1}) a
resolution of singularities $\overline X_n$ of an $n$-sheeted
cyclic covering, branched along $\bar H_{k}$, has the first Betti
number $b_1(\overline X_{k,n})=2k$. In addition, we show that for
any $k\in \mathbb N$, there is a Hurwitz curve $\bar H_{k}$
consisting of two irreducible components such that the resolution
of singularities $\overline X_{k,6}$ of a cyclic covering of the
projective plane, branched along $\bar H_{k}$, has the first Betti
number $b_1(\overline X_{k,6})=k$. Recall that $b_1(\overline
X_n)$ is always even if $\bar H$ is a plane algebraic curve, hence
$\bar H_k$ cannot be algebraic if $k$ is odd. Recall also that
Moishezon (\cite{Moi2}) proved the exist\-ence of an infinite
sequence $\bar H_i$ of irreducible cuspidal Hurwitz curves of
degree $54$ with exactly 378 cusps and 756 nodes which have
pairwise distinct braid monodromy type. In particular, they are
pairwise non-isotopic, and almost all of them are not isotopic to
an algebraic cuspidal curve.

\begin{cor} Let $\overline X_n$ be a resolution of singularities of a cyclic covering
of the projective plane branched along a Hurwitz curve $\bar H$
consisting of $k$ irreducible components and, maybe, along
$L_{\infty}$. If $n$ is divisible by $\deg \bar H$, then
$b_1(\overline X_n)=\deg \Delta (t)-k+1$.
\end{cor}
 \begin{cor} Let $\overline X_n$ be a resolution of singularities of a cyclic covering
of the projective plane of any degree $n$ branched along an
irreducible Hurwitz curve $\bar H$ and, maybe, along $L_{\infty}$.
If $\deg \bar H=p^k$, where $p$ is a prime number, then
$b_1(\overline X_n)=0$.
\end{cor}
\begin{cor} Let $\overline X_{p^k}$ be a resolution of singularities of a cyclic covering
of the projective plane of degree $p^k$ branched along any
irreducible Hurwitz curve $\bar H$ and, maybe, along $L_{\infty}$,
where $p$ is a prime number. Then $b_1(\overline X_n)=0$.
\end{cor}
Note that for any $k\in \mathbb N$, we show that there is a Hurwitz curve $\bar
H_k$ consisting of $k+1$ components and which is the branch curve
 of a 2-sheeted cyclic covering a resolution of singularities $\overline X_{k,2}$ which has
$b_1(\overline X_{k,2})=k$ (see Proposition \ref{b1G2}). In
particular, in our example the Hurwitz curve $\bar H_1$ has $\deg
\bar H_1=2^{10}$, the number of singular points of $\bar H_1$ is
equal to $2^{16}$, and all its singular points are of the form
$w^4-z^4=0$.\\

Recently, Auroux and Katzarkov (see \cite{Au}, \cite{Au-Ka})
proved the following theorem. Let $(\overline X,\omega)$ be a
compact symplectic 4-manifold with symplectic form $\omega$ with
class $[\omega]\in H^2(\overline X,\mathbb Z)$. Fix an
$\omega$-compatible almost complex structure $J$ and the corresponding
Riemannian metric $g$. Let $L$ be a line bundle on $\overline X$
whose first Chern class is $[\omega]$. Then, for $k>>0$, the line
bundle $L^{\otimes k}$ admits many approximately holomorphic
sections so that one can choose three of them which give an
approximately holomorphic generic covering $f_k:\overline X\to
\mathbb C\mathbb P^2$ of degree $N_k=k^2 \omega^2$ branched over a
cuspidal Hurwitz curve $\bar H_k$ (possibly, with negative nodes).

Any such covering $f_k:\overline X\to \mathbb C\mathbb P^2$ of
degree $N_k$, branched over a cuspidal Hurwitz curve $\bar H$,
determines a monodromy ${\mu}$, that is, an epimorphism
$\mu:\pi_1(\mathbb C^2\setminus H)\to \Sigma_{N_k}$ to the
symmetric group $\Sigma_{N_k}$ with additinal properties of genericity. On the
other hand, any homomorphism
$\mu:\pi_1(\mathbb C^2\setminus H)\to \Sigma_{N}$, such that
$\mu(\pi_1)$ acts transitively on a set consisting of $N$ elements,
defines an unramified covering $f:X\to \mathbb C^2\setminus H$ of
degree $N$. The covering $f$ can be extended to a covering
$\widetilde f:\widetilde X\to \mathbb C\mathbb P^2$ branched over
the Hurwitz curve $\bar H$ and, maybe, over $L_{\infty}$. In this paper we prove
(see Corollary \ref{sympl-cor}) that if $\widetilde X$ has arbitrary
analytic singularities (and if $\bar{H}$ has negative nodes, we assume that
the covering space is non--singular over them), then a resolution $\overline X$ of
singularities of $\widetilde X$ can be equipped with  a symplectic
structure.

The proofs of Theorems \ref{main}, \ref{Roots} and Corollary
\ref{cor} are given in section \ref{section5} and section
\ref{section6} is devoted to the proof of Theorem \ref{th}.

{\it Acknowledgement.}  The second author would like to express
his gratitude to University of Kaiserslautern for its hospitality
during the preparation of this paper.

\section{Representation of Hurwitz curves as sections of line
bundles} \label{sect1}

We begin with the following lemma.
\begin{lem} \label{definition}
The definitions of Hurwitz curves in $\mathbb C\mathbb P^2$ given
in \cite{Kh-Ku} and in the Introduction are equivalent.
\end{lem}
\proof Recall the definition of Hurwitz curves given in
\cite{Kh-Ku}. Let $F_1$ be a relatively minimal ruled rational
surface, $\mbox{pr}: F_1 \to \mathbb C\mathbb P^1$ the ruling, $R$
a fiber of $\mbox{pr}$, and $E_1$ the exceptional section,
$E_1^2=-1$. Identify $\mbox{pr}: F_1 \to \mathbb C\mathbb P^1$
with a linear projection $\mbox{pr}:\mathbb C\mathbb P^2\to
\mathbb C\mathbb P^1$  with center at a point $p\in \mathbb
C\mathbb P^2$ ($p$ is the blow down of $E_1$ to the point).

The image $\bar H=f(S)\subset F_1$ of a smooth map $f:S \to
F_1\setminus E_1$ of an oriented closed real surface $S$ is called
a Hurwitz curve {\rm (}with respect to {\rm $\mbox{pr}$)} of
degree $m$ if there is a finite subset $Z\subset\bar H$ such that:
\begin{itemize}
\item[(i)] $f$ is an embedding  of the surface $S\setminus
f^{-1}(Z)$ and for any $s\notin Z$, $\bar H$ and the fiber
$R_{\text{pr}(s)}$ of pr meet at $s$ transversely and with
positive intersection number; \item[(ii)] for each $s\in Z$ there
is a neighbourhood $U\subset F_1$ of $s$ such that $\bar H\cap U$
is a complex analytic curve, and the complex orientation of $\bar
H\cap U \setminus \{s\}$ coincides with the orientation
transported from $S$ by $f$; \item[(iii)]  the restriction of {\rm
$\mbox{pr}$} to $\bar H$ is a finite map of degree $m$.
\end{itemize}

To show that the definition in the Introduction implies the
definition in \cite{Kh-Ku}, let us perform several monoidal
transforms with centers at singular points of $\bar H$ (and at the
singularities of the proper transforms of $\bar H$) to resolve all
singular points of $\bar H$. Denote by $\sigma: \widetilde{\mathbb
C\mathbb P}^2\to \mathbb C\mathbb P^2$ the composition of these
monoidal transformations and by $S$ the proper transform of $\bar
H$. Then $S$ is a smooth real surface and $f=\sigma_{\mid S}$ is a
smooth map. To define an orientation on $S$, let us choose an
orientation at each non-critical point $p$ of $\text{pr}_{\mid
\bar H}$ so that the local intersection number of $\bar H$ and the
fibre $R$ passing through $p$ at the point $p$ will be equal to
$+1$. Obviously, these orientations are compatible for all
non-singular points of $\bar H$. Since near singular points this
orientation coincides with the orientation given by the complex
analytic structure (recall that $\bar H$ is complex analytic near
critical points of $\text{pr}_{\mid \bar H}$), this orientation
can be extended to the preimages of these critical points.

To show that the definition in \cite{Kh-Ku} implies the definition
in the Introduction, let us choose a fibre $R$ of $\text{pr}$ and
put $\mathbb C^2=F_1\setminus (R\cup E_1)$. Let $(u,v)$ be
coordinates in $\mathbb C^2$ such that the restriction of
$\text{pr}$ is given by $(u,v)\to u$. Let  $(u,v_1(u)),\dots,
(u,v_m(u))$ be the coordinates of the intersection points of $\bar
H$ and the fibre $R$ of $\text{pr}$ over a non-critical value $u$.
Consider
\begin{equation} \label{ad}
F(u,v)=\prod_{i=1}^m(v-v_i(u)).
\end{equation}
Obviously, the function $F(u,v)$, defined everywhere outside the
fibres over critical values, is smooth and can be extended to a
function on all $\mathbb C^2$ satisfying the properties of the
definition given in the Introduction. \edvo

Let a Hurwitz curve $\bar H_0$ of degree $m$ be given by equations
(\ref{function}). A smooth isotopy $h_t:\mathbb C\mathbb P^2\times
[0,1]\to \mathbb C\mathbb P^2\times [0,1]$ is called a $H$-{\it
isotopy} if for each $t\in [0,1]$ the image $\bar H_t=h_t(\bar
H_0)$ is a Hurwitz curve given by equations
$$v_i^m
+\sum_{j=0}^{m-1}c_{j,i}(u_i,t)v_i^j=0,\qquad i=1,2,$$
$c_{j,i}(u_i,0)=c_{j,i}(u_i)$ for all $i,j$. (Note that in the
definition of $H$-isotopy given in \cite{Kh-Ku}, it is assumed
that the number of critical values of $\bar H_t$ does not depend
on $t$.) It is easy to see that if $\bar H_0$ and $\bar H_1$ are
$H$ isotopic and a line $L_{\infty}$ is generic with respect to
both Hurwitz curves $\bar H_0$ and $\bar H_1$, then $\mathbb C^2
\setminus H_0$ and $\mathbb C_2\setminus H_1$ are diffeomorphic.

Denote by $p_{\infty}=\mathbb C\mathbb P^2\setminus (\mathbb
C_1^2\cup\mathbb C_2^2)$ the center of the projection $\text{pr}$.
In what follows we will assume that the fibre of $\text{pr}$ over
$u_2=0$ is generic with respect to $\bar H_0$. Denote it by
$L_{\infty}$. Obviously, there is a smooth $H$-isotopy $h_t$
identical outside a small neighbourhood $U$ of $L_{\infty}$ such
that the function $F_2(u_2,v_2,1)$ defining $\bar H_1=h_1(\bar
H_0)$ in $\mathbb C_2^2$ coincides with the function $v_2^m-1$ at
the points $(u_2,v_2)$ with $\mid u_2\mid <\varepsilon$ for some
$\varepsilon
>0$. In what follows we will assume that
\begin{equation} \label{assump1}
\begin{array}{l} L_{\infty} \, \, is\, \,  given\, \,  by\, \, u_2=0\, \,
and\, \, \bar H \, \, is\, \,  given\, \,  by\, \, the\, \,
equation \\ v_2^m-1=0 \, \, in\, \, a\, \,
neighbourhood\, \, of \, \,  L_{\infty}. 
\end{array}
\end{equation}

Let $u_{1,j}$ be a critical value of a Hurwitz curve $\bar H_0$ of
degree $m$ given in $\mathbb C_1^2$ by equation $F_1(u_1,v_1)=0$,
that is, the number of the different roots of the equation
\begin{equation} \label{aaa} F_1(u_{1,j},v_1,0)=0 \end{equation}
is less than $m$, and let $v_{1,j_0}$ be a root of equation
(\ref{aaa}) of multiplicity one. It is obvious that there is a
smooth $H$-isotopy $h_t$ identical outside a small neighbourhood
$U=\{ \mid u_1-u_{1,j}\mid <\varepsilon\}$ such that the function
$F_1(u_1,v_1,1)$ defining $\bar H_1=h_1(\bar H_0)$ in $\mathbb
C_1^2$ is such that $v_1=v_{1,j_0}$ is a root of the equation
$F_1(u_{1},v_1,1)=0 $ for all $u_{1}$ such that $\mid
u_{1}-u_{1,j}\mid <\varepsilon_1$ for some positive
$\varepsilon_1<\varepsilon$. Therefore, in what follows we can
(and we will) assume that if $u_{1,j}$ is a critical value of
$\bar H$, then
\begin{equation} \label{assump2}
 \begin{array}{l}there\, \, is\, \, an\, \varepsilon >0\, \, such\, \, that\, \,
F_1(u_1,v_1)\, \, defining \, \, \bar H \,  \,  is \\  analytic\,
\, at \, \, (u_{1},v_1)\, \, for \, \, \mid u_{1}-u_{1,j}\mid
<\varepsilon \, \, and \, \, all \, \, v_1.
\end{array}
\end{equation}

Let us consider a line bundle $p:\mathcal L(k)\to \mathbb C\mathbb
P^2$ associated with the sheaf $\mathcal O_{\mathbb C\mathbb
P^2}(k)$. Recall its definition. The projective plane $\mathbb
C\mathbb P^2$ with homogeneous coordinates $(z_0:z_1:z_2)$ is
covered by three charts $\mathbb C_i^2$, $i=1,2,3$, isomorphic to
$\mathbb C^2$, with coordinates $(u_i,v_i)$, $u_1=z_1/z_0$,
$v_1=z_2/z_0$, $u_2=z_0/z_1$, $v_2=z_2/z_1$, $u_3=z_0/z_2$,
$v_3=z_1/z_2$. The bundle $\mathcal L(k)$ is covered by three
charts $W_i=\mathbb C_i^2\times \mathbb C_i^1$ with the third
coordinate $w_i$, $w_1=w_2/u_2^k$, $w_1=w_3/u_3^k$,
$w_2=w_3/v_3^k$, and the restriction of $p_{\mid W_i}$ coincides
with the projection to the first factor.

\begin{lem}\label{section}
The functions $w_i=F_i(u_i,v_i)$, $i=1,2$, defining $\bar H$
define a smooth section $s$ of $\mathcal L(m)$ over $\mathbb
C\mathbb P^2\setminus \{ p_{\infty}\}$.
\end{lem}
\proof In $\mathbb C_1^2\cap\mathbb C_2^2$ we have
$$\begin{array}{ll}
F_1(u_1,v_1)= & \displaystyle v_1^m
+\sum_{j=0}^{m-1}c_{j,1}(u_1)v_1^j=(\frac{v_2}{u_2})^m
+\sum_{j=0}^{m-1}c_{j,1}(\frac{1}{u_2})(\frac{v_2}{u_2})^j= \\
\displaystyle & \displaystyle (\frac{1}{u_2})^m(v_2^m
+\sum_{j=0}^{m-1}c_{j,1}(\frac{1}{u_2})u_2^{m-j}v_2^j).
\end{array}
$$
The functions $$F_2(u_2,v_2)=v_2^m
+\sum_{j=0}^{m-1}c_{j,2}(u_2)v_2^j$$ and $$v_2^m
+\sum_{j=0}^{m-1}c_{j,1}(\frac{1}{u_2})u_2^{m-j}v_2^j$$ coincide
with each other, since they are smooth and for almost all (except
a finite number of) values $u_{2,0}$ of $u_2$ the polynomials
$$v_2^m
+\sum_{j=0}^{m-1}c_{j,2}(u_{2,0})v_2^j$$ and $$v_2^m
+\sum_{j=0}^{m-1}c_{j,1}(\frac{1}{u_{2,0}})u_{2,0}^{m-j}v_2^j$$
have the same sets of roots. \edvo
\begin{lem} \label{homotopy} Let $f_0:S^3\to \mathbb C^*$ be a smooth
function on $S^3=\{ (u,v)\in \mathbb C^2\mid \, u\overline
u+v\overline v=\varepsilon^2\}$, $0<\varepsilon<<1$, such that
$f_0$ coincides with the function $v^m-1$ in a neighbourhood
$U\subset S^3$ of the circle $u=0$. Then there is a smooth
function $F:S^3\times [0,1]\to \mathbb C^*$ such that
\begin{itemize}
\item[($i$)] $F(u,v,0)=f_0(u,v)$;
\item[($ii$)] $F(u,v,t)=v^m-1$
for $(u,v)\in U$ and  $t\in [0,1]$;
\item[($iii$)]
$F(u,v,1)=v^m-1$.
\end{itemize}
\end{lem}
\proof  Since $S^3$ is simply connected, there is a lift
$\widetilde f_0:S^3\to \widetilde{\mathbb C}^*$ of the function
$f_0$ such that $f_0=e\circ \widetilde f_0$, where $\widetilde{
\mathbb C}^*$ is the complex plane $\mathbb C$ with complex
coordinate $x$ and $e:\widetilde{ \mathbb C}^*\to \mathbb C^*$ is
the universal covering given by $y=e^x$. Without loss of
generality one can assume that $\widetilde f_0(0,\varepsilon)=\ln
(1-\varepsilon^m)+\pi i$. Denote by $f_1:S^3\to \mathbb C^*$ the
function $v^m-1$ and by $\widetilde f_1:S^3\to \widetilde{\mathbb
C}^*$ its lift such that $\widetilde f_1(0,\varepsilon)=\ln
(1-\varepsilon^m)+\pi i$. Therefore we have $\widetilde f_{0\mid
U}\equiv \widetilde f_{1\mid U}$.

Consider a function $\widetilde F:S^3\times [0,1]\to \mathbb C$
given by $$x=t\widetilde f_1(u,v)+(1-t)\widetilde f_0(u,v).$$
Obviously, the function $F=e\circ \widetilde F$ has the  desired
properties. \edvo
\begin{lem} \label{sect2} There is a real number $\varepsilon_1$,
$0<\varepsilon_1 <<1$, and a smooth section $\overline s_m$ of
$\mathcal L(m)$ over $\mathbb C\mathbb P^2$ such that
\begin{itemize} \item[($i$)]
$\bar H\subset \mathbb C\mathbb P^2\setminus B(\varepsilon_1)$,
where $B(\varepsilon_1)=\{ u_3\overline u_3+v_3\overline v_3\leq
\varepsilon_1^2\}$ is a ball in $\mathbb C\mathbb P^2$ with center
at $p_{\infty}$; \item[($ii$)] over $\mathbb C\mathbb P^2\setminus
B(\varepsilon_1)$ the section $\overline s_m$ coincides with the
section $s$ from Lemma {\rm \ref{section}}; \item[($iii$)]
$\overline s_m$ is complex analytic in a neighbourhood of the line
$L_{\infty}$.
\end{itemize}
\end{lem}
\proof By the definition of Hurwitz curves, there is a ball
$B(\varepsilon_1)=\{ u_3\overline u_3+v_3\overline v_3\leq
\varepsilon_1^2\}$ for some positive $\varepsilon_1$ such that
$\bar H\subset \mathbb C\mathbb P^2\setminus B(\varepsilon_1)$.

The line bundle $\mathcal L_m$ is trivial over $B(\varepsilon_1)$.
Therefore the restriction of the section $s$ from Lemma
\ref{section} to $\partial B(\varepsilon_1)=S^3$ defines a
function $f_0:S^3\to \mathbb C^*$. Denote by $F:S^3\times [0,1]\to
\mathbb C^*$ a function existence of which is proved in Lemma
\ref{homotopy} (in the notations of Lemma \ref{homotopy}, we put
$u=u_3$, $v=v_3$ and $\varepsilon=\varepsilon_1$). Let us choose
$\varepsilon_2<\varepsilon_1$ and a smooth monotone function
$r:[\varepsilon_2,\varepsilon_1]\to [0,1]$ such that
$r(\varepsilon_1)=0$ and $r(\varepsilon_2)=1$. Put
$h:B(\varepsilon_1,\varepsilon_2)=B(\varepsilon_1)\setminus
B(\varepsilon_2)\to
\partial B(\varepsilon_1)$ to be the map given by $$h(u_3,v_3)=
(\frac{\varepsilon_1u_3}{\sqrt{u_3\overline u_3+v_3\overline
v_3}}, \frac{\varepsilon_1v_3}{\sqrt{u_3\overline u_3+v_3\overline
v_3}}),$$
$\widetilde F(u_3,v_3,t)=h^*(F)$, and
$$\overline
F(u_3,v_3)=(\frac{\sqrt{u_3\overline u_3+v_3\overline
v_3}}{\varepsilon_1})^m(\widetilde F(u_3,v_3,r(\sqrt{u_3\overline
u_3+v_3\overline v_3}))+1)-1.$$ Then the section $\widetilde s$,
defined as follows,
\[ \widetilde s(p)=\left\{
\begin{array}{ll} s(p) \qquad \qquad & \text{for} \,\, p\in \mathbb C\mathbb P^2\setminus B(\varepsilon_1)  \\
\overline F(u_3,v_3)\qquad & \text{for} \,\, p\in
B(\varepsilon_1)\setminus
B(\varepsilon_2) \\
v_3^m-1  \qquad \qquad & \text{for} \,\, p\in  B(\varepsilon_2),
\end{array} \right.
\]
satisfies all conditions of Lemma \ref{sect2}, except that,
possibly, it is not smooth, but only continuous at the points from
$B=(\partial B(\varepsilon_1)\cup
\partial B(\varepsilon_2))\setminus U$, where $U$ is a neighbourhood
of $L_{\infty}$. By standard theorems from analysis, there is a
smooth section $\overline s_m$ close enough to $\widetilde s$
which coincides with $\widetilde s$ outside a small enough
neighbourhood $V$ of $B$ such that $\overline V\cap (\bar H\cup
L_{\infty})=\emptyset$, where $\overline V$ is the closure of $V$.
\edvo

\section{Symplectic varieties with analytic singularities 
} \label{section2}

Let $Y$ be a projective complex manifold, $\dim_{\mathbb C} Y=n$
and $\omega$ a K\"{a}hler form on $Y$, $[\omega ]\in H^2(Y,\mathbb
Z)$. Consider $(Y,\omega)$ as a symplectic manifold,
$\dim_{\mathbb R} Y=2n$. A closed subvariety $X$ of $Y$ is called
a {\it symplectic variety with analytic singularities} if there
are open subsets $U_0\subset U\subset Y$ such that the closure
$\overline U_0$ in $Y$ is a subset of $U$, $X\cap U$ is a complex
analytic subset in $U$ and $X\setminus \overline U_0$ is a smooth
symplectic submanifold. Denote by $\text{Sing}\, X$ the set of
points of $X$ in which $X$ is not smooth. Then $\text{Sing}\, X$
is a projective algebraic subvariety of $Y$.

\begin{lem} \label{sympl-lem} Let $X$ be a symplectic variety with analytic
singularities in a projective complex manifold $Y$ with
K\"{a}hler form  $\omega$.  Let $Z\subset \text{Sing}\, X$ be a
nonsingular projective subvariety of $Y$, $\sigma :\overline Y\to
Y$ the monoidal transformation of $Y$ with center in $Z$, and
$\overline X$ be the proper transform of $X$. Then there is a
K\"{a}hler form $\overline{\omega}$ on $\overline Y$ such that
$\overline X$ is a symplectic subvariety with analytic
singularities of $(\overline Y,\overline \omega)$.
\end{lem}
\proof The manifold $\overline Y$ is projective algebraic.
Consider an embedding $i:\overline Y\hookrightarrow \mathbb
C\mathbb P^N$ to some projective space and denote by $\varphi
=i\circ \sigma^{-1}$ a rational map from $Y$ to $\mathbb C\mathbb
P^N$. Let $\Gamma\subset Y\times\mathbb C\mathbb P^N$ be the
closure of the graph of $\varphi$ and $p_i$, $i=1,2$, be the
projections from $Y\times\mathbb C\mathbb P^N$ to the factors. The
morphisms $i$ and $\sigma$ define the morphism $\sigma \times i
:\overline Y\to \Gamma\subset Y\times\mathbb C\mathbb P^N$. Since
the composition $p_2\circ (\sigma \times i):\overline Y\to
\overline Y$ is an isomorphism, $p_{2\mid \Gamma}:\Gamma \to
\overline Y$ is an isomorphism too. Moreover, if we identify
$\overline Y$ with $\Gamma$ by means of $p_{2\mid \Gamma}$, then
$p_{1\mid \Gamma}$ coincides with $\sigma$.

Denote by $\Omega=\Omega_N$,
$$\Omega_N=\frac{i}{(\sum_{j=0}^N\overline z_jz_j)^2}\sum_{k=0}^N\sum_{j\neq k}
(\overline z_jz_jdz_k\wedge d\overline z_k -\overline
z_jz_kdz_j\wedge d\overline z_k),$$ the Fubini -- Studi symplectic
form on $\mathbb C\mathbb P^N$. It is a K\"{a}hler integer form.
Consider the restriction of a form
$\omega_{\varepsilon}=p^*_1(\omega)+\varepsilon p^*_2(\Omega)$ to
$\Gamma$. It is a K\"{a}hler form on $\Gamma$ for each
$\varepsilon
>0$.

Choose open neighbourhoods $V_0\subset V\subset Y$ of $\text{Sing}\, X$
such that $V\cap X$ is an analytic subvariety and the closure
$\overline V_0$ of $V_0$ in $Y$ is a subset of $V$. Denote by
$X_0=X\setminus V_0$. It is  compact and $\sigma^{-1}_{\mid
X_0}:X_0\to \overline X_0=\sigma^{-1}(X_0)$ is an isomorphism.
Therefore $\overline X_0$ is compact.

Evidently, the restriction of $\omega _{\varepsilon}$ to
$\Gamma\cap p_1^{-1}(V)$ is a symplectic form at each nonsingular
point of $\Gamma\cap p_1^{-1}(V)$ for all $\varepsilon >0$, since
$\Gamma\cap p_1^{-1}(V)$ is an analytic set in $p_1^{-1}(V)$.
Since the restriction of $\omega $ to $X_0$ is a symplectic form
at each point of $X_0$ and $\overline X_0$ is  compact, we can
choose $\varepsilon$ to be small enough, so that the restriction
of $\omega_{\varepsilon}=p^*_1(\omega)+\varepsilon p^*_2(\Omega)$
to $\overline X_0=\Gamma\cap p_1^{-1}(X_0)$ is a symplectic form
at each point of $\overline X_0$. If we take $\varepsilon =\frac{m}{n}$ rational, then $n\omega_{\varepsilon}$ is an
integer form. \edvo

\section{Symplectness of coverings of the projective plane branched
along Hurwitz curves} \label{section3}

In this section we use the notations and assumptions of section
\ref{sect1}.

Let $\bar H$ be a Hurwitz curve, possibly with {\it negative
nodes}, that is, in a neighbourhood $U$ of each critical point
$p$, it is either given by an analytic equation or $\bar H\cap U$
consists of two smooth branches meeting transversely at $p$ with
intersection number $-1$, and each branch of $\bar H\cap U$ meets
the fiber $\text{pr} ^{-1}(\text{pr}(p))$ transversely at $p$ with
intersection number equals $+1$.

We fix a point $p\in \mathbb C\mathbb P^2\setminus (\bar H
\cup L_{\infty})$. Consider the fundamental group
$\pi_1=\pi_1(\mathbb C^2\setminus H, p)$ of the complement of the
affine Hurwitz curve $H=(\mathbb C\mathbb P^2\setminus
L_{\infty})\cap \bar H$.  Let us choose a point $x\in \bar
H\setminus \text{Sing}\, \bar H$  and consider a line $L \subset
\mathbb C^2$ meeting $H$ transversely at $x$. Let $\gamma \subset
L$ be a circle of small radius with center at $x$. The choice of
an orientation on $\mathbb C^2$ defines an orientation on $\gamma
$. Let $\Gamma $ be a loop consisting of a path $l$ in $\mathbb
C^2\setminus H$, joining $p$ with a point $q\in \gamma$, the loop
$\gamma $ (with positive direction) starting and ending at $q$,
and a return path to $p$ along $l$ in the opposite direction. Such a
loop $\Gamma$ (and the corresponding element in $\pi _1$) is
called a {\it geometric generator} (with center at $x$) of the
fundamental group $\pi _1=\pi _1(\mathbb C^2\setminus H,\, p)$. It
is well known that $\pi _1$ is generated by geometric
generators.

For each critical point $s_i$ of $H$ let us choose a neighbourhood
$U_i\subset \mathbb C^2$ such that $H\cap U_i$ is given (in local
coordinates in $U_i$) by an analytic equation or, if $s_i$ is a
negative node, then it consists of two smooth branches meeting
transversely at $p$. Note that if $s_i$ is a negative node, then
$\pi _1(U_i\setminus H,p_i)$ is isomorphic to $\mathbb Z\oplus
\mathbb Z$ and generated by two commuting geometric generators.

Choose smooth paths $\gamma _i$ lying in $\mathbb C^2\setminus H$
and connecting the points $p_i$ with $p$. This choice defines
homomorphisms $\psi _i:\pi _1(U_i\setminus H,p_i)\to \pi _1$. We
call  $\psi _i(\pi _1(U_i\setminus H,p_i))=G_i$  {\it the local
fundamental group} of the singular point $s_i$. The local
fundamental groups are defined uniquely up to conjugation in
$\pi_1$.

Consider a homomorphism $\mu :\pi_1\to \Sigma_N$ from the
fundamental group $\pi_1=\pi_1(\mathbb C^2\setminus H)$ of the
complement of the affine Hurwitz curve $H=(\mathbb C\mathbb
P^2\setminus L_{\infty})\cap \bar H$ to the symmetric group
$\mathfrak S_N$ such that its image $Im\,\mu$ acts transitively on
a set consisting of $N$ elements.

Let $s_i$ be a negative node of $H$. As it was mentioned above,
the local fundamental group $G_i$ is generated by two commuting
geometric generators, say $\Gamma _{i,1}$ and $\Gamma_{i,2}$.
Denote by $N_{i,j}=\{1\leq n\leq N\mid \mu (\Gamma_{i,j})(n)\neq
n\}$. We say that $\mu$ is {\it good at the negative node $s_i$}
if $N_{i,1}\cap N_{i,2}=\emptyset$. The homomorphism $\mu $ is
called a {\it monodromy of degree $N$} if it is good at all
negative nodes.

The homomorphism $\mu$ defines an unramified covering $f=f_{\mu}:
Y\to \mathbb C^2\setminus H$  of degree $N$. This covering can be
extended to a finite ramified covering $\widetilde f :\widetilde Y
\to \mathbb C\mathbb P^2$ branched along $\bar H$ and, maybe,
along $L_{\infty}$.

To describe this extension, consider a geometric generator $\Gamma
$ with center at $x\in H\setminus \text{Crit}\, H$, where
$\text{Crit}\, H$ is the set of critical points of $H$. The image
$\mu (\Gamma)$ in $\Sigma_N$ is a product of cyclic permutations
$\sigma_{1},\dots, \sigma_{n_x}$ (it is possible that the orders
of some $\sigma_l$ are equal to one). Let
$\sigma_l=(n_{1,l},\dots,n_{r_{l},l})$ be a permutation of order
$r_{l}$ , where $1\leq n_{j,l}\leq N$, then the number of the
preimages $\widetilde f^{-1}(x)$ is equal to $n_x$ and each point
$y$ from $\widetilde f^{-1}(x)$ corresponds to a cyclic
permutation $\sigma_l$. Near the point $y_l$ corresponding to a
cyclic permutation $\sigma _l$, the covering $\widetilde f$ is a
cyclic covering of degree $r_l$ branched along $H$ and it is
locally isomorphic to a subvariety of $\mathbb C^3$ given by
$w^{r_l}= v-v_j(u)$, where $v-v_j(u)=0$ is a local equation of
$\bar H$ at the point $x$ (see (\ref{ad})). These local
isomorphisms define on $\widetilde Y$ a structure of a smooth
manifold at each point $y$ lying over $H\setminus \text{Crit}\,
H$.

Let $s_i\in \text{Crit} \, H$ be a negative node. As it was
mentioned above, the local fundamental group $G_i$ is generated by
two geometric generators $\Gamma_{i,1}$ and $\Gamma_{i,2}$. The
images $\mu (\Gamma_{i,j})$ in $\Sigma_N$ are  products of cyclic
permutations $\sigma_{1,i,j},\dots, \sigma_{k_{i,j},i,j}$. Let
$\sigma_{e,i,j}= (n_{1,l,i,j},\dots, n_{r_{l,i,j},l,i,j})$ be  a
permutation of order $r_{l,i,j}$. Put $N_{i,1,2}=\{ 1\leq n\leq
N\mid \mu(\Gamma_{i,1})(n)=n\, \, \text{and}\, \,
\mu(\Gamma_{i,2})(n)=n\}$. Since $\mu$ is a monodromy, the set
$\widetilde f^{-1}(s_i)$ is in one-to-one correspondence with the
union of the cyclic permutations $\sigma_{l,i,j}$, $j=1,2$, of
orders bigger than one and the set $N_{i,1,2}$. Moreover, if $y\in
\widetilde f^{-1}(s_i)$ corresponds to an element from
$N_{i,1,2}$, then $\widetilde f$ is an isomorphism of a
neighbourhood $V$ of $y$ and its image $\widetilde f(V)$, and if
$y\in \widetilde f^{-1}(s_i)$ corresponds to a cyclic permutation
$\sigma_{l,i,j}$ of order bigger than one, then the restriction of
$\widetilde f$ to a neighbourhood of $y$ is a cyclic covering of a
neighbourhood of $s_i$ of degree $r_{l,i,j}$ branched along the
$j$-th  branch of the negative node. It is locally isomorphic to a
subvariety of $\mathbb C^3$ given by $w^{r_{l,i,j}}= v-v_j(u)$,
where $v-v_j(u)=0$ is a local equation of the $j$-th branch of
$\bar H$ at the point $s_i$. These local isomorphisms define on
$\widetilde Y$ a structure of smooth manifold at each point $y\in
\widetilde f^{-1}(s_i)$.

Let $x=s_i\in \text{Crit}_{\text{analytic}} H$, that is, $s_i$ is
a critical point of $H$ which is not a negative node. Then over a
small neighbourhood $U$ of $s_i$ in which $H$ is given by an
analytic equation, the preimage $\overline f^{-1}(U)$ is the
disjoint union of $n_{s_i}$ open neighbourhoods being in
one-to-one correspondence with the orbits of the action of $\mu
(G_i)$ on the set consisting of $N$ elements. By the theorem of
Grauert, Remmert and Stein (for a proof see \cite{St}) we know
that, over a neighbourhood $U$ of $s_i$ the variety $\widetilde Y$
can be equipped with the structure of a two dimensional complex
analytic variety.

By assumption, in a neighbourhood $U$ of $L_{\infty}$ (where
$U=\mathbb C\mathbb P^2\setminus B(R)$ and $B(R)\subset \mathbb
C^2_1$ is a ball of big radius $R$) the curve $\bar H$ coincides
with the algebraic curve $\bar C\subset \mathbb C\mathbb P^2$ of
degree $m$ given in $U$ by the equation $v_2^m-1=0$. If we choose the
base point $p$ to lie in $U$, we can consider  a Zariski -- van
Kampen presentations of $\pi_1$ and of the fundamental group
$\widetilde \pi_1=\pi_1(\mathbb C\mathbb P^2\setminus (\bar C\cup
L_{\infty}),p)$ having the same sets of generators and it is easy
to see that these presentations define an epimorphism $e
:\widetilde\pi_1\to\pi_1$. The composition $\mu\circ e$ defines a
ramified covering $\widetilde g:\widetilde Z\to \mathbb C\mathbb
P^2$ branched along $\bar C$ and, maybe, along $L_{\infty}$. It is
easy to see that the coverings $\widetilde f$ and $\widetilde g$
are isomorphic over $U$. Therefore the variety $\widetilde
f^{-1}(U)$ can be identified with $\widetilde g^{-1}(U)$ by means
of an isomorphism $h:\widetilde f^{-1}(U)\to \widetilde g^{-1}(U)$
so that $\widetilde f^{-1}(U)$ can also be considered as a complex
analytic variety. Let $i: \widetilde Z\hookrightarrow \mathbb
C\mathbb P^{m_{\infty}}$ be an embedding such that $\widetilde g$
is defined by the projection $(z'_0:z'_1:\dots :z'_{m_{\infty}})\to
(z'_0:z'_1:z'_2)$. Put $z_j=h^*(z'_j)$, $j=3,\dots,m_{\infty}$ and
\begin{equation} \label{coor}
w'_{j,\infty}=\frac{z_j}{z_0},\qquad j=3,\dots,m_{\infty}.
\end{equation}
\begin{thm} \label{cov-sym}
Let $\bar H$ be a Hurwitz curve with negative nodes, $\mu
:\pi_1(\mathbb C^2\setminus H)\to \Sigma_N$ a monodromy, and
$\widetilde f:Y\to \mathbb C\mathbb P^2$ the covering
associated with the monodromy $\mu$. Then $\widetilde Y$ can be
embedded into some projective space $\mathbb C\mathbb P^M$ as a
symplectic subvariety with analytic singularities.
\end{thm}

\proof Below, for each point $p\in \mathbb C^2_1$, we denote  by $V_p\subset
U_p$ small balls $V_p=\{ \mid u_1-u_1(p)\mid^2+\mid
v_1-v_1(p)\mid^2< \delta_1\}$ and $U_p=\{ \mid
u_1-u_1(p)\mid^2+\mid v_1-v_1(p)\mid^2< \delta_2\}$ of radii
$\delta_1$ and $\delta_2$, $0<\delta_1<\delta_2<<1$, and by
$\rho_{p
}:\mathbb C\mathbb P^2\to \mathbb R$ we denote a smooth
non negative function such that $\rho_{p
\mid V_p}\equiv 1$ and $\rho_{p \mid \mathbb C\mathbb P^2\setminus
U_p}\equiv 0$.

To construct the desired embedding, let us choose two open coverings $\{
U_i\}$ and $\{ V_i\}$ of $\mathbb C\mathbb P^2$ as follows.

For each point $s_i\in \text{Crit}_{\text{analytic}}\bar H$ which
is not a negative node, let us choose small neighbourhoods
$V'_{s_i}\subset V_{s_i}\subset U_{s_i}\subset U'_{s_i}$  such
that
\begin{itemize} \item[($c_1$)] the curve $U'_{s_i}\cap \bar H$ is
analytic in $U'_{s_i}$; \item[($c_2$)] the preimage $\widetilde
f^{-1}(U'_{s_i})$ splits into the disjoint union of neighbourhoods
 of the points
$y_{i,j}\in \widetilde f^{-1}(s_i)$; \item[($c_3$)] the radius of
the ball $V'_{s_i}$ (resp. $U_{s_i}$) is strictly less than the
radius of $V_{s_i}$ (resp. $U'_{s_i}$).
\end{itemize}

Let $V'_{\infty}\subset V_{\infty}\subset U_{\infty}\subset
U'_{\infty}$ be open neighbourhoods of $L_{\infty}$ such that
\begin{itemize}
\item[($c_4$)] in $U'_{\infty}$, the Hurwitz curve $\bar H$
coincides with the curve $\bar C$ given in $U'_{\infty}$ by
equation $v_2^m-1=0$; \item[($c_5$)]  $U'_{\infty}\cap
V'_{s_i}=\emptyset$ for all neighbourhoods $V'_{s_i}$ of singular
points $s_i$ chosen above; \item[($c_6$)] the closure $\overline
V'_{\infty}\subset V_{\infty}$ and the closure $\overline
U_{\infty}\subset U'_{\infty}$. \end{itemize}

Let $\rho_{\infty}:\mathbb C\mathbb P^2\to \mathbb R$ be  a smooth
non negative function such that $\rho_{\infty \mid
V_{\infty}}\equiv 1$ and $\rho_{\infty \mid \mathbb C\mathbb
P^2\setminus U_{\infty}}\equiv 0$.

Let us add the neighbourhoods $U_{\infty}$ and $V_{\infty}$ to the
sets $\{ U_{s_i}\}$ and $\{ V_{s_i}\}$ chosen above.

Next, for each point $p\in \mathbb C^2_1\setminus ((\cup
V_{s_i})\cup V_{\infty})$, we can find open neighbourhoods $V_{p}\subset
U_{p}$ of $p$ such that
\begin{itemize}  \item[($c_7$)]  $U_{p}\cap V'_{s_i}=\emptyset$ and $U_{p}\cap
V'_{\infty}=\emptyset$ for the neighbourhoods $V'_{s_i}$ and
$V'_{\infty}$ chosen above; \item[($c_8$)] the preimage
$\widetilde f^{-1}(U_{p})$ splits into the disjoint union of
neighbourhoods of the points $y_{j}\in \widetilde f^{-1}(p)$.
\item[($c_9$)] if $p\in U'_{s_i}$ (resp. $p\in U'_\infty)$ for
some neighbourhood $U'_{s_i}$ (resp. $U'_\infty)$ chosen above,
then $U_p\subset
  U'_{s_i}$ (resp. $U_p \subset U'_\infty$);
\item[($c_{10}$)] if $p\not\in (\cup_{s_i\in
\text{Crit}_{\text{analytic}}\bar H}U'_{s_i})\cup U'_{\infty}$,
then $U_p\cap ((\cup_{s_i\in \text{Crit}_{\text{analytic}}\bar
H}U_{s_i})\cup U_{\infty})=\emptyset$.
\end{itemize}

Let us add the neighbourhoods $U_{p}$ and $V_{p}$ to the sets $\{
U_{s_i}\}$ and $\{ V_{s_i}\}$ chosen above (here, one of
$s_i=\infty$). As a result, we obtain two open coverings $\mathcal
V= \{ V_p\}$ and $\mathcal U=\{ U_p\}$ of $\mathbb C\mathbb P^2$.

For $U'=U'_{\infty}$ let $(w_{3,\infty},\dots
,w_{m_{\infty},\infty})$ be functions in $\widetilde Y\setminus
\widetilde f^{-1}(L_{\infty})$ defined by $w_{j,\infty}=\widetilde
f^*(\rho_{\infty})w'_{j}$ for $j=3,\dots,m_{\infty}$, where the
functions $w'_j$ were defined in (\ref{coor}).

For each $U'_{s_i}$, where $s_i\in \text{Crit}_{\text{analytic}}
\bar H$, there are complex analytic functions $w'_{1,s_i},\dots
,w'_{m_{s_i},s_i}$ in $\widetilde U'_{s_i}=\widetilde
f^{-1}(U'_{s_i})$ such that these functions together with
$\widetilde f^*(u_1)$ and $\widetilde f^*(v_1)$ give an analytic
embedding of $\widetilde U'_{s_i}$ to $\mathbb C^{m_{s_i}+2}$.
Denote by $w_{j,s_i}=\widetilde f^*(\rho_{s_i})w'_{j,s_i}$, $1\leq
j\leq m_{s_i}$, functions in $\widetilde Y$.

By construction of the open coverings, the preimage $\widetilde
f^{-1}(V_p)=\coprod \widetilde V_{p,j}$ (resp. $\widetilde
f^{-1}(U_p)=\coprod \widetilde U_{p,j}$) splits into disjoint
union of $m_p=n_p$ connected neighbourhoods $\widetilde V_{p,j}$
(resp. $\widetilde U_{p,j}$), $j=1,\dots, n_p$.  If $p\in \bar H$,
then a neighbourhood $\widetilde U_{p,j}$ is isomorphic to a
subvariety of $\mathbb C^3$ given in coordinates
$(u_1,v_1,w'_{j,p})$ by
$$(w'_{j,p}-w^0_{j,p})^{r_{j,p}}=v_1-v_{1,p}(u_1),$$ where
$v_1-v_{1,p}(u_1)=0$ is the equation of $\bar H$ in $U_p$. Extend
the functions $w'_{j,p}$ by putting $w'_{j,p\mid \widetilde
U_{p,l}}\equiv 0$ for $l\neq j$ and choose constants $w^0_{j,p}$,
$j=1,\dots,m_p$, so that the functions $(u_1,v_1,w'_{1,p},\dots,
w'_{m_p,p})$ define a smooth embedding of $\widetilde U_p$ to
$\mathbb C^{m_p+2}$. Denote by $w_{j,p}=\widetilde
f^*(\rho_p)w'_{j,p}$, $j=1,\dots,m_p$, functions in $\widetilde
Y$.

Note that if $y\in \widetilde f^{-1}(\bar H)\cap \widetilde
V_{p,j}\cap \widetilde U_{q,j}$, where $p,q\neq \infty$ and
$p,q\not\in \text{Crit}_{\text{analytic}}\bar H$, then by the
definition of the functions $w_{j,p}$ and $w_{j,q}$, we have
$r_{j,p}=r_{j,q}=r_j$ and there is a $r_j$-th root $\zeta_{p,q}$
of unity such that
\begin{equation} \label{connection1}\begin{array}{ll}
w_{j,q} & =\rho_q(u_1,v_1)(\zeta_{p,q}(w_{j,p}-w^0_{j,p})+w^0_{j,p}) \\
w_{j',q} & \equiv 0 \qquad \qquad \qquad \qquad \qquad \qquad
\text{for}\, \, j'\neq j
\end{array}
\end{equation}
in a neighbourhood of $y$.

Similarly, if $y\in \widetilde f^{-1}(\bar H)\cap \widetilde
V_{p,j}\cap \widetilde U_{q}$ (resp. if $y\in \widetilde
f^{-1}(\bar H)\cap \widetilde V_{q}\cap \widetilde U_{p,j}$),
where $q= \infty$ or $q\in \text{Sing}_{\text{analytic}}\bar H$,
then by the definition of the functions $w_{j,p}$ and $w_{i,q}$
and by properties $(c_1)$, $(c_4)$, $(c_9)$, and $(c_{10})$, we
have
\begin{equation} \label{connection2}
(w_{j,p}-w_{j,p}^0)^{r_j}=v_1-F(u_1)
\end{equation}
and
\begin{equation} \label{connection3}
w_{i,q}=\rho_q(u_1,v_1)h_{i}(u_1,v_1,w_{j,p}) \qquad \text{for}\,
\, 1\leq i\leq m_q
\end{equation}
(resp. $w_{j,p}
=\rho_p(u_1,v_1)h_{j}(u_1,v_1,w_{1,q},w_{2,q},w_{3,q},\dots,w_{m_q,q})$,
here $w_{1,q}$ and $w_{2,q}$ are constants if $q=\infty$) in a
neighbourhood of $y$, where $F$ and all $h_i$ (resp. $h_j$) are
analytic functions, and $v_1-F(u_1)=0$ is an analytic equation of
some branch of $\bar H$.

If $p\not\in \bar H\cup L_{\infty}$, then we can assume that
$\widetilde f$ defines an isomorphism of the neighbourhoods
$\widetilde U_{p,j}$ and $U_p$ for $j=1,\dots, m_p=N$. Choose $N$
different constants $w^0_{j,p}$ and define functions
$w_{j,p}=\widetilde f^*(\rho_p)w'_{j,p}$, $j=1,\dots,m_p$, in
$\widetilde Y$, where the function $w'_{j,p}$ is defined in
$\widetilde U_p$ as follows
\[
 w'_{j,p}(q)\equiv \left\{ \begin{array}{ll}
w^0_{j,p}\qquad & \text{if}\, \,q\in \widetilde U_{p,j} \\
0\qquad & \text{if}\, \,q\not\in \widetilde U_{p,j}.
\end{array}
\right.
\]

Choose a finite covering $\widetilde{\mathcal V}_0=\{ \widetilde
V_{p_i,j}\mid 1\leq i\leq k, \, \, 1\leq j\leq m_{p_i}\}\cup
\{\widetilde V_{\infty}\}$ of $\widetilde Y$, set
$$M=m_{\infty}+\sum_{j=1}^km_{p_j}, $$
$w_3=w_{3,\infty},\dots,w_{m_{\infty}}=w_{m_{\infty},\infty}$, and
enumerate the set of functions
$$\{ w_{j,p_i}\mid 1\leq i\leq k, \, \, 1\leq j\leq m_{p_i}\}
$$ by the numbers $m_{\infty}+1,\dots, M$.

Consider a linear projection $p:\mathbb C\mathbb P^M\to \mathbb
C\mathbb P^2$ given by $$(z_0:z_1:z_2:\dots :z_M)\to
(z_0:z_1:z_2).$$ The base locus of $p$ is a projective space
$P\simeq \mathbb C\mathbb P^{M-3}$ given by equations
$z_0=z_1=z_2=0$. The restriction of $p$ to $\mathcal L=\mathbb
C\mathbb P^M\setminus P$ defines on $\mathcal L$ the  structure of a
vector bundle over $\mathbb C\mathbb P^2$ the zero section of
which is given by $z_3=\dots =z_M=0$. Over the charts $\mathbb
C_i^2$ with coordinates $(u_i,v_i)$, $i=1,2,3$, the bundle
$\mathcal L$ is trivial, it is isomorphic to $\mathbb C_i^M\simeq
\mathbb C_i^2\times \mathbb C^{M-2}_i$. In particular,
$(z_3/z_0,\dots,z_M/z_0)$ are coordinates in $\mathbb C_1^{M-2}$.

Over $\mathbb C_1^2$, consider a map $\alpha' :\widetilde
f^{-1}(\mathbb C_1^2)\to \mathbb C_1^M$ given by
$$\alpha' (y)=(\widetilde f^*(u_1)(y),\widetilde f^*(v_1)(y),w_{3}(y),\, .\, .\, .\, ,w_{M}(y)).$$

Since for $p_i\neq \infty$ each $U_{p_i}$ is a subset of $\mathbb
C_1^2$ and all functions $w_{j,p_i}\equiv 0$ at the points lying
over the complement of $U_{p_i}$, the map $\alpha'$ can be
extended to a map $\alpha: \widetilde Y\to \mathcal L$ as follows.
Over the neighbourhood $\widetilde f^{-1}(V'_{\infty})$ the map
$\alpha \text{ equals }i\circ h$, where $h$ was defined above and
$i$ is a linear embedding of $\mathbb C\mathbb P^{m_{\infty}}$ to
$\mathbb C\mathbb P^M$ given by $$i((z_0: \, .\, .\, .\,
:z_{m_{\infty}}))=(z_0: \, .\, .\, .\, :z_{m_{\infty}}:0:\, . \,
.\, .\, :0).$$ It is easy to see that $\alpha$ is an embedding
such that $\alpha (\widetilde V'_{s_i})$ and $\alpha (\widetilde
V'_{\infty})$ are analytic subsets of $\mathcal L$ for $\widetilde
V'_{\infty}=\widetilde f^{-1}(V'_{\infty})$ and for all
neighbourhoods $\widetilde V'_{s_i}=\widetilde f^{-1}(V'_{s_i})$,
$s_i\in \text{Crit}_{\text{analytic}}\bar H$.

Denote by $\Omega$ the restriction of the Fubini -- Studi form
$\Omega_M$ to $\mathcal L$. In the chart $\mathbb C_1^M$, it has
the form
\begin{equation} \label{Omega}
\Omega=\frac{\displaystyle i\sum_{k=1}^M (dw_k\wedge d\overline
w_k +\displaystyle \sum_{j\neq k} (\overline w_jw_jdw_k\wedge
d\overline w_k -\overline w_jw_kdw_j\wedge d\overline
w_k))}{\displaystyle (1+\sum_{j=1}^M\overline w_jw_j)^2},
\end{equation}
where $w_k=\frac{z_k}{z_0}$ and $w_1=u_1$, $w_2=v_1$.

Denote by the same symbol $\overline \varepsilon =(\varepsilon_1,
\varepsilon_2)$ a collection of two positive numbers and the
linear transformation $\overline \varepsilon :\mathbb C\mathbb
P^M\to \mathbb C\mathbb P^M$ given by
$$(z_0:z_1:z_2:z_3: \, .\, .\, .\, :z_M)\to
(z_0:z_1:\varepsilon_1z_2: \varepsilon_2z_3:\, .\,. \, .\,
:\varepsilon_{2}z_M).$$ Denote by $\omega_{\overline \varepsilon}$
the restriction of the form $\Omega$ to $\widetilde Y_{\overline
\varepsilon}=(\overline\varepsilon\circ\alpha)(\widetilde Y)$. Let
us show that there exist a positive constant $c_1$ and a positive
function $c_2(t)$, $t\in (0,c_1]$, such that $\widetilde
Y_{\overline \varepsilon}$ is a symplectic in $\mathcal L$ for all
$\overline \varepsilon =(\varepsilon _1,\varepsilon_2)$ with
$\varepsilon_1\leq c_1$, $\varepsilon _2\leq c_2(\varepsilon_1)$.
Note that if $\varepsilon _1$ is small enough, then the image of
$\bar H$ under the map $(z_0:z_1:z_2)\to
(z_0:z_1:\varepsilon_1z_2)$ becomes symplectic.

 For each
$\overline \varepsilon$, the form $\omega_{\overline \varepsilon}$
is a symplectic form at the points from the neighbourhoods
$(\overline \varepsilon \circ \alpha )(\widetilde V'_{s_i})$,
where $s_i$ is an analytic singular point of $\bar H$, and from
$(\overline \varepsilon \circ \alpha )(\widetilde V'_{\infty})$,
since $\widetilde Y_{\overline \varepsilon}$ is an analytic
subvariety of $\mathcal L$ at these points.

For each point $y$ not belonging to the ramification locus of
$\widetilde f$, the variety $\widetilde Y$ is given locally at
$\alpha (y)$ by equations $w_j=F_j(u_1,v_1)$, $j=3,\dots, M$,
where $F_j(u_1,v_1)$ are smooth functions in a neighbourhood of
$\widetilde f(y)$. Therefore the variety $\widetilde Y_{\overline
\varepsilon}$ is given locally at $(\overline \varepsilon \circ
\alpha )(y)$ by equations
$w_j=\varepsilon_{2}F_j(u_1,\frac{v_1}{\varepsilon_1})$,
$j=3,\dots, M$. It is easy to see that for each fixed
$\varepsilon_1$ and for $\varepsilon _{2}$ being small enough, the
form $\omega_{\overline\varepsilon}$ is symplectic at $(\overline
\varepsilon \circ \alpha )(y)$, since the variety $\widetilde
Y_{\overline \varepsilon}$ is very close to the algebraic variety
given by $w_j=0$, $j=3,\dots, M$, which is symplectic.

Consider a point $y\in \widetilde f^{-1}(\bar H)$ belonging to the
ramification locus of $\widetilde f$ and such that $\widetilde
f(y)\not\in \widetilde U_p$ for $p\in
\text{Crit}_{\text{analytic}}\bar H$ and for $p={\infty}$. By
(\ref{connection1}), renumbering the coordinates $w_3,\dots,
w_{M}$, we can assume that $\widetilde Y_{\overline \varepsilon}$
is given in a neighbourhood of $\alpha (y)$ by equations
\begin{equation} \label{equation3}
\begin{array}{rll}(w_3-w_{3,0})^r
& = v_1-F(u_1), &
\\ w_j & =
\rho_j(u_1,v_1) h_j(w_3), & j\geq 4,
\end{array}
\end{equation}
where $r\geq 2$, $\rho_j$ are smooth functions,
$h_j=\zeta_{3,j}(w_3-w_{3,0})+w_{j,0}$ are analytic functions,
$v_1-F(u_1)=0$ is the equation of a branch of $\bar H$ at the
point $\widetilde f(y)$, and the point $\alpha
(y)=(u_{1,0},F(u_{1,0}),w_{3,0},\dots,w_{M,0})$. Then the variety
$\widetilde Y_{\overline \varepsilon}$ is given by equations
\begin{equation} \label{equation}
\begin{array}{rll}\varepsilon_1(w_3-\varepsilon_2w_{3,0})^r
& = \varepsilon_2^r(v_1-\varepsilon_1F(u_1)), &
\\ w_j & = \varepsilon_2\rho_j(u_1,\frac{v_1}{\varepsilon_1})
h_j(\frac{w_3}{\varepsilon_2}), & j\geq 4.
\end{array}
\end{equation}
in a neighbourhood of $(\overline \varepsilon\circ\alpha
)(y)=(u_{1,0},\varepsilon_1F(u_{1,0}),\varepsilon_2w_{3,0},\dots,\varepsilon_2w_{M,0})$.

Denote by $A_1=\frac{\partial F}{\partial u_1}(\alpha (y))$,
$A_2=\frac{\partial F}{\partial \overline u_1}(\alpha (y))$,
$B_j=\frac{\partial h_j}{\partial w_3}(\alpha (y))$,
$C_{j,1}=\frac{\partial \rho _j}{\partial u_1}(\alpha (y))$,
$C_{j,2}=\frac{\partial \rho _j}{\partial \overline u_1}(\alpha
(y))$, $D_{j,1}=\frac{\partial \rho _j}{\partial v_1}(\alpha
(y))$, $D_{j,2}=\frac{\rho _j}{\partial \overline v_1}(\alpha
(y))$, $\rho_{j,0}=\rho_j(\alpha(y))$, $h_{j,0}=h_j(w_{3,0})$,
$j=4,\dots, M$. It follows from (\ref{equation}) that at the point
$(\overline \varepsilon \circ \alpha )(y)$, we have

\begin{equation} \label{dif}
\begin{array}{lll} dv_1 & = & \varepsilon_1(A_1du_1+A_2d\overline u_1) \\
d\overline v_1 & = & \varepsilon_1(\overline A_{2}du_1+\overline
A_{1}d\overline u_1) \end{array} \end{equation} and
\begin{equation} \label{dif2}
\begin{array}{lll}
dw_j & = & \rho_{j,0}B_jdw_3+ \varepsilon_2h_{j,0}(C_{j,1}du_1+
C_{j,2}d\overline u_1
+D_{j,1}\frac{dv_1}{\varepsilon_1}+D_{j,2}\frac{d\overline
v_1}{\varepsilon_1})
 \\
d\overline w_j & = & \rho_{j,0}\overline B_jd\overline w_3+
\varepsilon_2\overline h_{j,0}(\overline C_{j,2}d\overline u_1+
\overline C_{j,1}du_1 +\overline D_{j,2}\frac{d
v_1}{\varepsilon_1}+\overline D_{j,1}\frac{d\overline
v_1}{\varepsilon_1}),
\end{array}\end{equation}
$j=4,\dots,M$.

If we substitute   (\ref{dif}) to (\ref{dif2}), we obtain that

\begin{equation} \label{dif3}
\begin{array}{lll}
dw_j & = & \rho_{j,0}B_jdw_3+ \varepsilon_2\nu _j
 \\
d\overline w_j & = & \rho_{j,0}\overline B_jd\overline w_3+
\varepsilon_2\overline \nu_j,
\end{array}\end{equation}
where the forms $\nu_j$ and $\overline \nu_j$ do not depend on
$\overline \varepsilon$ for $j=4,\dots,M$.

It follows from (\ref{dif}) and (\ref{dif3}) that for each very
small $\overline \varepsilon $ the tangent space of $\widetilde
Y_{\overline \varepsilon}$  at the point $(\overline \varepsilon
\circ \alpha )(y)$ is very close to the tangent space at the point
$(\overline \varepsilon \circ \alpha )(y)$ of a linear algebraic
variety $Z$ given by $v_1=\varepsilon_1F(u_{1,0})$,
$w_j-w_{j,0}=\rho_{j,0}B_j(w_3-w_{3,0})$, $j=4,\dots,M$. Therefore
for each very small $\overline \varepsilon $ the form
$\omega_{\overline \varepsilon}$ is symplectic at $(\overline
\varepsilon \circ \alpha )(y)$. By continuity, it is
symplectic in some neighbourhood of  $(\overline \varepsilon \circ
\alpha )(y)$.

Consider a point $y\in \widetilde f^{-1}(\bar H)$ belonging to the
ramification locus of $\widetilde f$ and such that $\widetilde
f(y)\in U_p$ for some $p\in \text{Crit}_{\text{analytic}}\bar H$
or $p={\infty}$. By (\ref{connection2}) and (\ref{connection3}),
renumbering the coordinates $w_3,\dots, w_{M}$, we can assume that
there is some $n$, $3\leq n\leq M$, such that $\widetilde
Y_{\overline \varepsilon}$ is given in a neighbourhood of $\alpha
(y)$ by equations
\begin{equation} \label{equation5}
\begin{array}{rll} h_j(u_1,v_1,w_3,\dots,w_{n}) & =0, &  j=3,\dots, n
\\
\rho_j(u_1,v_1)h_j(u_1,v_1,w_3,\dots,w_{n}) & =w_j , &
j=n+1\dots,M,
\end{array}
\end{equation}
where $\rho_j$ are smooth functions and $h_j$ are analytic
functions at the point $\widetilde f(y)$. Let
$(u_{1,0},v_{1,0},w_{3,0},\dots,w_{M,0})$ be the coordinates of
the point $\alpha (y)$. Then the variety $\widetilde Y_{\overline
\varepsilon}$ is given in a neighbourhood of $(\overline
\varepsilon\circ\alpha
)(y)=(u_{1,0},\varepsilon_1v_{1,0},\varepsilon_2w_{3,0},\dots,\varepsilon_2w_{M,0})$
by equations
\begin{equation} \label{equation6}
\begin{array}{rll} h_j(u_1,\frac{v_1}{\varepsilon_1},\frac{w_3}{\varepsilon_2},\dots,\frac{w_{n}}{\varepsilon_2})
& =0, &  j=3,\dots, n
\\
\varepsilon_2\rho_j(u_1,\frac{v_1}{\varepsilon_1})
h_j(u_1,\frac{v_1}{\varepsilon_1},\frac{w_3}{\varepsilon_2},\dots,\frac{w_{n}}{\varepsilon_2})
& =w_j , & j=n+1\dots,M.
\end{array}
\end{equation}

Set $A_{j,l}=\frac{\partial h_j}{\partial w_l}(\alpha
(y)
)\,\, \text{for} \, j\leq j\leq M, \, 1\leq l\leq n,$ (here
$w_1=u_1$ and $w_2=v_1$) and by $B_j=\rho_j(u_{1,0},v_{1,0})$ for
$n+1\leq j\leq M$.

It follows from (\ref{equation6}) that for each fixed positive
$\varepsilon_1$ and for very small $\varepsilon_2$ the tangent
space of $\widetilde Y_{\overline \varepsilon}$  at the point
$(\overline \varepsilon \circ \alpha )(y)$ is very close to the
tangent space at the point $(\overline \varepsilon \circ \alpha
)(y)$ of a linear algebraic variety $Z$ given by
$$\begin{array}{l}
\displaystyle \sum_{l=3}^nA_{j,l}(w_l-\varepsilon_2w_{l,0})
=-A_{j,1}(u_1-u_{1,0})-\frac{A_{j,2}}{\varepsilon_1}(v_1-\varepsilon_1v_{1,0}),
\, \,  3\leq j\leq n
\\
\displaystyle B_j\sum_{l=3}^nA_{j,l}(w_l-\varepsilon_2w_{l,0})
=w_j-\varepsilon_2w_{j,0}, \qquad n+1\leq j\leq M.
\end{array}
$$
Therefore for fixed $\varepsilon_1$ and for each very small
$\varepsilon _2$, the form $\omega_{\overline \varepsilon}$ is
symplectic at $(\overline \varepsilon \circ \alpha )(y)$. By
continuity, it is symplectic in some neighbourhood of
$(\overline \varepsilon \circ \alpha )(y)$.

To complete the proof, it suffices to apply the compactness of
$\widetilde Y$.
 \edvo
\medskip

By Lemma \ref{sympl-lem} and by Hironaka's Theorem on resolution
of singularities, we have
\begin{cor} \label{sympl-cor}
Let $\widetilde f :\widetilde Y\to \mathbb C\mathbb P^2$ be a
finite covering branched along a Hurwitz curve $\bar H$ (possibly,
with negative nodes) and, maybe, along $L_{\infty}$, and
associated with a monodromy $\mu :\pi_1(\mathbb C^2\setminus H)\to
\Sigma_N$. Then there exist collections of positive integers
$(M_1,\dots,M_k)$ and $(n_1,\dots,n_k)$ such that a resolution
$\overline Y$ of singularities of $\widetilde Y$ can be embedded
as a symplectic submanifold to $(\mathbb C\mathbb P^{M_1}\times
\dots \times \mathbb C\mathbb P^{M_k},\Omega_{n_1,\dots,n_k})$,
where
$\Omega_{n_1,\dots,n_k}=n_1p_1^*(\Omega_{M_1})+\dots+n_kp_k^*(\Omega_{M_k})$
and $\Omega_{M_j}$ is the Fubini -- Studi symplectic form on
$\mathbb C\mathbb P^{M_j}$.
\end{cor}

\begin{thm} In the notations of the proof of Theorem {\rm \ref{cov-sym}},
let $\widetilde Y$ be a smooth manifold. Then the symplectic
structure, constructed in the proof of Theorem {\rm \ref{cov-sym}}
and given by the symplectic form $\omega_{\overline \varepsilon}$,
does neither depend on $\overline \varepsilon$ if the coordinates of
$\overline \varepsilon$ are small enough nor on the choice the
coverings $\mathcal U$ and $\mathcal V$, and on the choice of
the functions $w_{i,j}$.

Moreover, if $i:\widetilde Y\hookrightarrow \mathbb C\mathbb P^N$
is an algebraic embedding (in the case when $\bar H$ is an
algebraic curve) and $\widetilde f=p\circ i$, where $p:\mathbb
C\mathbb P^N\to \mathbb C\mathbb P^2$ is a linear projection, then
$(\widetilde Y, \omega_{\overline \varepsilon})$ and $(\widetilde
Y, i^*(\Omega_N))$ are symplectomorphic for $\overline
\varepsilon$ being small enough, where $\Omega_N$ is the Fubini --
Studi form on $\mathbb C\mathbb P^N$ and $\omega_{\overline
\varepsilon}$ is the form constructed in the proof of Theorem {\rm
\ref{cov-sym}}.
\end{thm}
\proof Let $\omega_{\overline \varepsilon}=(\overline
\varepsilon\circ \alpha)^*(\Omega)$, where the embedding $\alpha
:\widetilde Y\to \mathcal L$ was constructed in the proof of
Theorem \ref{cov-sym}. Note that for each $\overline \varepsilon$,
$0<\varepsilon_1 \leq c_1$ and $0<\varepsilon_2\leq
c_2(\varepsilon_1)$, the class $[\omega_{\overline
\varepsilon}]\in H^2(\widetilde Y, \mathbb Z)$ is dual to the
class $[\widetilde f^{-1}(L)]\in H_2(\widetilde Y, \mathbb Z)$,
where $L$ is a line in $\mathbb C\mathbb P^2$. Therefore, by
Moser's stability theorem for symplectic structures (see
\cite{McD-S}, Theorem 3.17), the forms $\omega_{\overline
\varepsilon}$ define the same symplectic structure if
$0<\varepsilon_1 \leq c_1$ and $0<\varepsilon_2\leq
c_2(\varepsilon_1)$.

The symplectic structure on $\widetilde Y$, defined by the forms
$\omega_{\overline \varepsilon}$ does not depend on the choice of
the coverings $\{ U_i\}$ and $\{ V_i\}$, and the choice of the
functions $w_{i,j}$ defining the embedding $\alpha$. Indeed, let
two collections $\{ w'_{i,j}\}$ and $\{ w^{\prime\prime}_{i,j}\}$
of functions define two embeddings $\alpha ':\widetilde Y\to
\mathcal L'\subset \mathbb C\mathbb P^{M'}$ and $\alpha
^{\prime\prime}:\widetilde Y\to \mathcal L^{\prime\prime}\subset
\mathbb C\mathbb P^{M^{\prime\prime}}$, and put
$\omega_{\overline\varepsilon}'=(\overline \varepsilon
\circ\alpha')^*(\Omega ')$ and $\omega_{\overline
\varepsilon}^{\prime\prime}=(\overline \varepsilon
\circ\alpha^{\prime\prime})^*(\Omega ^{\prime\prime})$, where
$\Omega ^{\prime}$ and $\Omega ^{\prime\prime}$ are the Fubini --
Studi symplectic forms on $\mathbb C\mathbb P^{M'}$ and $\mathbb
C\mathbb P^{M^{\prime\prime}}$ respectively. Then we have an
embedding $\alpha'\times \alpha^{\prime\prime} :\widetilde Y\to
\mathcal L'\times_{\mathbb C\mathbb P^2}\mathcal
L^{\prime\prime}$. Note that the form $\Omega_t=t (p')^*(\Omega
')+(1-t)(p^{\prime\prime})^*(\Omega^{\prime\prime})$ is a
K\"{a}hler form for each $t\in [0,1]$. Since the segment $[0,1]$
is compact, applying a similar calculation as in the proof of
Theorem \ref{cov-sym}, one can show that there exists
$\overline\varepsilon
=\overline\varepsilon'=\overline\varepsilon^{\prime\prime}$,
$0<\varepsilon_1\leq \min (c_1',c_1^{\prime\prime})$,
$0<\varepsilon_2\leq \min
(c'_2(\varepsilon_1),c_2^{\prime\prime}(\varepsilon_1))$, such
that $\omega_{t,\overline\varepsilon}=((\overline\varepsilon \circ
(\alpha'\times \alpha^{\prime\prime}))^*(\Omega_t)$ is a
symplectic form on $\widetilde Y$ for all $t\in [0,1]$.  On the
other hand,
$\omega_{0,\overline\varepsilon}=\omega'_{\overline\varepsilon}$
and
$\omega_{1,\overline\varepsilon}=\omega^{\prime\prime}_{\overline\varepsilon}$,
and the forms $\omega_{t,\overline\varepsilon}$ belong to the same
cohomology class. Therefore, by Moser's stability theorem for
symplectic structures, the forms $\omega_{t,\overline\varepsilon}$
define the same symplectic structure on $\widetilde Y$.

In the case of an algebraic embedding $i:\widetilde
Y\hookrightarrow \mathbb C\mathbb P^N$ such that $\widetilde
f=p\circ i$, where $p:\mathbb C\mathbb P^N\to \mathbb C\mathbb
P^2$ is a linear projection, denote by $\alpha ':\widetilde Y\to
\mathcal L'\subset \mathbb C\mathbb P^{M'}$ any embedding
constructed in the proof of Theorem \ref{cov-sym} and put $\alpha
^{\prime\prime}=i$. Applying the same arguments as above, it is
easy to see that $(\widetilde Y, \omega_{\overline \varepsilon})$
and $(\widetilde Y, i^*(\Omega_N))$ are symplectomorphic for
$\overline \varepsilon$ being small enough, since the symplectic
manifolds $(\widetilde Y,\omega_{\overline
\varepsilon}^{\prime\prime})$ are symplectomorphic for all
positive collections $\overline \varepsilon$, where
$\omega_{\overline \varepsilon}^{\prime\prime}=(\overline
\varepsilon \circ\alpha^{\prime\prime})^*(\Omega _N)$. \edvo

\section{Embeddings of cyclic coverings of the plane into rational projective 3-folds }
\label{section4}

In this section, we use the notations and assumptions of section
\ref{sect1}.

Let $\bar H$ be a Hurwitz curve of degree $m$. Consider the
infinite cyclic covering $f =f_{\infty }:X_{\infty } \to X'={\Bbb
C}^{2}\, \setminus \, H$ corresponding to the epimorphism $\nu
:\pi_1(\mathbb C^2\setminus H)\to \mathbb F_1$. The covering
$f_{\infty}$ can be factorized through the cyclic covering
$f_n:X'_n\to \mathbb C^2\setminus H$ associated with the
epimorphism $\text{mod}_n\circ \nu$, $f_{\infty}=g_n\circ f_n$.

In this section, we will show that the covering $f_n$ can be
extended to a smooth map ${\overline f}_n:{\overline X}_n\to
\mathbb C\mathbb P^2$ branched along $\bar H$ and, maybe, along
$L_{\infty}$ (if $n$ is not a divisor of $\deg \bar H$, then
${\overline f}_n$ is branched along $L_{\infty}$), where
$\overline X_n$ is a real smooth 4-fold.

\begin{thm} \label{sympl-thm} A resolution of singularities $\overline X_n$ of a cyclic covering
of $\mathbb C\mathbb P^2$ of degree $n$, branched along a Hurwitz
curve $\bar H$ and, maybe, along $L_{\infty}$, can be embedded into
some rational projective 3-fold (equipped with integer K\"{a}hler
symplectic structure) as a symplectic submanifold.
\end{thm}

\proof Since $\mathbb C^2\setminus H_1$ and $\mathbb C^2\setminus
H_2$ are diffeomorphic for $H$-isotopic Hurwitz curves $\bar H_1$
and $\bar H_2$, we can assume that $\bar H$ satisfies conditions
(\ref{assump1}) and (\ref{assump2}).

By Lemma \ref{sympl-lem} and by Hironaka's Theorem on resolution
of singularities, it suffices to show that for some extension
$\widetilde f_n:\widetilde X_n\to \mathbb C\mathbb P^2$ of
$f'_n:X'_n\to X'$, the variety $\widetilde X_n$ can be imbedded to
some rational projective 3-fold (equipped with integer K\"{a}hler
symplectic structure) as a symplectic subvariety with analytic
singularities.

To show this, denote by $d$ the smallest non-negative integer for
which $m+d$ is divisible by $n$. Put $m+d=kn$ and consider the  line
bundle $p:\mathcal L(k)\to \mathbb C\mathbb P^2$ (see section
\ref{sect1}) associated with the sheaf $\mathcal O_{\mathbb
C\mathbb P^2}(k)$. By Lemma \ref{sect2}, the Hurwitz curve $\bar
H$ coincides with the zero locus of
 a smooth section $\overline s_m$ of
$\mathcal L(m)$ over $\mathbb C\mathbb P^2$ such that $\overline
s_m$ is analytic in a neighbourhood $U$ of $L_{\infty}$ and in
neighbourhoods of all critical points of $\bar H$.

 Denote by $\overline s_d$ the
section of $\mathcal L(d)$ defined over $\mathbb C_2^2$ by
$w_2=u_2^d$. The product
\begin{equation} \label{sectionmd} \overline
s_{m+d}=\overline s_m\overline s_d \end{equation} is a section of
$\mathcal L(m+d)$, where $\overline s_m$ is a section of $\mathcal
L(m)$ satisfying all conditions from Lemma \ref{sect2}.

Define $\alpha :\widetilde X_n\hookrightarrow  \mathcal L(k)$ by the
equation
\begin{equation} \label{X_n}
w_i^n=\overline s_{m+d}(u_i,v_i)
\end{equation}
and put $\widetilde f_n=p_{\mid X_n}$, where $p:\mathcal L(k)\to
\mathbb C\mathbb P^2$ is the morphism defining on $\mathcal L(k)$
the structure of the line bundle.
 In
particular, $\widetilde X_n$ is given by the equation
$$w_1^n=F_1(u_1,v_1)$$
in $\mathbb C_1^3$ and it is given by
$$ w^n_2=u_2^dF_2(u_2,v_2)$$
in $\mathbb C_2^3$.

Obviously, over $\mathbb C\mathbb P^2\setminus (\bar H\cup
L_{\infty})$ the covering $\widetilde f_n$ is an unramified
$n$-sheeted cyclic covering. Next, all singular points of the
variety $\widetilde X_n$ lie over singular points of $\bar H$ and,
maybe, over $L_{\infty}$. Moreover, by construction of the section
$\overline s_{m+d}$, the set $\text{Sing}\, \widetilde X_n$ is
complex analytic in some  neighbourhood $U\subset \mathcal
L(k)$.

The line bundle $\mathcal L(k)$ is a quasi-projective variety and
it can be compactified to a projective three dimensional rational
manifold $\overline{\mathcal L}(k)$ by adding a section "at infinity".

The variety $\overline{\mathcal L}(k)$ has many different
embeddings to projective spaces, since its Picard group
$\text{Pic}(\overline{\mathcal L}(k))\simeq \mathbb Z\oplus
\mathbb Z$. We can choose one of these embeddings, for example,
the following one.

In the neighbourhood $\mathbb C^3_1$ with coordinates
$(u_1,v_1,w_1)$ consider monomials $u_1^{a_1}v_1^{a_2}w_1^{a_3}$,
$0\leq a_1+a_2+ka_3\leq k+1$, the number of which is equal to
$\frac{(k+2)(k+3)}{2}+3$. Put $N=\frac{(k+2)(k+3)}{2}+2$ and
consider a rational map $h:\overline{\mathcal L}(k)\to \mathbb
C\mathbb P^N$ given in $\mathbb C_1^3$ by $z_{\overline
a}=u_1^{a_1}v_1^{a_2}w_1^{a_3}$, where $\overline a=(a_1,a_2,a_3)$
are triples of integers and $z_{\overline a}$ are homogeneous
coordinates in $\mathbb C\mathbb P^N$. It is easily to check  that
$h$ is an embedding.

Consider the Fubini - Studi form $\Omega _N$ on $\mathbb C\mathbb
P^N$ and denote by $\Omega =h^*(\Omega_N)$ its pull back.

As in the proof of Theorem \ref{cov-sym}, denote by the same
symbol $\overline \varepsilon =(\varepsilon_1, \varepsilon_2)$ a
collection of two positive numbers and automorphism of $\overline
{\mathcal L}(k)$ given in $\mathbb C_1^3$ by $(u_1,v_1,w_1)\to
(u_1,\varepsilon_1v_1,\varepsilon_2w_1)$.

The calculations (we omit them) similar to the calculations done
in the proof of Theorem \ref{cov-sym} show that there exist a
positive constant $c_1$ and a positive function $c_2(t)$ such that
$\widetilde X_{\overline \varepsilon}=(\overline\varepsilon\circ
\alpha)(X_n)$ is a sympectic subvariety of $\mathcal L$ with
analytic singularities for all $\overline \varepsilon
=(\varepsilon _1,\varepsilon_2)$ with $\varepsilon_1\leq c_1$,
$\varepsilon _2\leq c_2(\varepsilon_1)$.
 \edvo

\section{Alexander polynomials
of Hurwitz $C$-groups}\label{section5}

Let $\bar H$ be a Hurwitz (resp. topological Hurwitz) curve of
degree $m$. Since any Zariski -- van Kampen presentation of
$\pi_1(\mathbb C^2\setminus H)$ is a $C$-presentation of a Hurwitz
$C$-group of degree $m$, Theorems \ref{main} and \ref{Roots} are
corollaries of the following Theorems \ref{Main} and \ref{roots}.

\begin{thm} \label{Main} Let $G\in \mathcal H$ be a Hurwitz $C$-group of degree $m$
and $\Delta (t)$ its Alexander polynomial. Then
\begin{itemize}
\item[($i$)] $\Delta (t)\in \mathbb Z[t]$; \item[($ii$)]
$\Delta(0)=\pm 1$; \item[($iii$)] the roots of $\Delta (t)$  are
$m$-th roots of unity; \item[($iv$)] the rank
of the free part of $N^{\prime}/N^{\prime\prime}$ is equal to $\deg
\Delta (t)$; \item[($v$)]  the action of $h_{\mathbb C}$ on
$N/N'\otimes \mathbb C$ is semisimple.\end{itemize}
\end{thm}

\proof Consider the exact sequence of groups
$$ 1\to N\to G \buildrel{\nu}\over\longrightarrow \mathbb F_1\to 1. $$
This exact sequence induces an automorphism $\widetilde h\in
\text{Aut}\, N$ (an action of the $C$-generator $x\in \mathbb F_1$
on $N$) given by $\widetilde h(n)= \widetilde x^{-1}n\widetilde x$
for $n\in N$, where $\widetilde x$ is one of the $C$-generators of
$G$. Evidently, the automorphism $\widetilde h$ is defined
uniquely up to inner automorphisms of the group $N$, therefore
$\widetilde h$ defines an automorphism $h\in \text{Aut} \, N/N'$.

In \cite{Ku.O}, it was proved that for a Hurwitz $C$-group $G$,
the group $N$ is finitely presented. Therefore $N/N'$ is a
finitely generated abelian group. Let $N/N'=T\oplus F$ be a
decomposition into the direct sum of the torsion subgroup $T$ and
a free abelian group $F$. Note that $T$ is a finite group and $F$
is finitely generated. The automorphism $h$ of $N/N'$ induces an
automorphism of $T$ and therefore an automorphism $\overline h$ of
$F\simeq (N/N')/T$. If one chooses a free basis of the $\mathbb
Z$-module $F$ over $\mathbb Z$, then this automorphism will be
given by a matrix $H$ with integer coefficients. Since the
automorphism $h_{\mathbb C}$ of $N/N'\otimes \mathbb C$ can be
given by the same matrix $H$, the polynomial $\Delta (t)=\det
(H-t\text{Id})\in \mathbb Z[t]$. And since $\overline h\in
\text{Aut}\, F$, then $\det \widetilde H=\pm 1$ and therefore,
$\Delta(0)=\pm 1$.

Let us show that $\widetilde h^m$ is an inner automorphism of $N$.
Indeed, since $G$ is a Hurwitz group of degree $m$, it is
generated by $C$-generators $x_1,\dots , x_m$ such that the
product $x_1\dots x_m$ belongs to the center of $G$. Therefore the
element $\widetilde x^m=\widetilde n\cdot x_1\dots x_m$, where
$\widetilde n\in N$, induces on $N$ an inner automorphism. Thus
the induced automorphisms $h^m$ of $N/N'$ and $h_{\mathbb C}^m$ of
$N/N'\otimes \mathbb C$ are trivial, that is, $h_{\mathbb
C}^m=\text{Id}$.

Since $h_{\mathbb C}^m=\text{Id}$, all roots of the polynomial
$\Delta(t)=\det (h_{\mathbb C}-t\text{Id})$ are $m$-th roots of
unity, and the action of $h_{\mathbb C}$ on $N/N'\otimes \mathbb
C$ is semisimple. \edvo

\begin{thm} \label{roots}
Let $G$ be an irreducible Hurwitz $C$-group, then
\begin{itemize}
\item[($i$)] $\Delta (1)=1$; \item[($ii$)] $\deg \Delta (t)$ is an
even number; \item[($iii$)] $\Delta (t)$ is a reciprocal
polynomial.
\end{itemize}
\end{thm}

\proof  It follows from Lemma 6 in \cite{Ku3} that $\Delta (1)=\pm
1$. Let us show that $t=-1$ also is not a root of $\Delta (t)$.
Indeed, if $t=-1$ is a root of $\Delta (t)$, then $\Delta
(t)=(t+1)P(t)$, where $P(t)$ is a polynomial with integer
coefficients. Therefore $2P(1)=\pm 1$. But this is impossible, since
$P(1)$ is an integer.

By Theorem \ref{Main}, all roots of $\Delta (t)\in \mathbb Z[t]$
are roots of unity. They are non-real, since $t=\pm 1$ are not
roots of $\Delta(t)$.  Thus, $\deg \Delta (t)$ is an even number
and ($ii$) is proved.

It is well known that if $\lambda$ is a primitive $k$-th root of
unity, $k>2$, and a polynomial $P(t)\in \mathbb Z[t]$ has
$\lambda$ as one of its roots, then all primitive $k$-th roots of
unity are roots of $P(t)$. In particular, $\lambda^{-1}$ is also a
root of $P(t)$ and ($iii$) is proved.

Since $\Delta (t)=\det (h_{\mathbb C}-t\text{Id})$ and $\deg
\Delta (t)$ is an even number, then $\Delta (t)= t^{\deg \Delta
(t)}+\dots $ . Let $\displaystyle \Delta (t)=\prod_i\Phi_{k_i}(t)$
be a factorization as the product of $k_i$-th cyclotomic
polynomials. Now statement ($i$) follows from well-known Lemma
\ref{LL}, since $\Delta (1)=\pm 1$.
\begin{lem} \label{LL}
Let $\Phi_{k}(t)$  be a  $k$-th cyclotomic polynomial, $k>1$. Then
\[ \Phi_{k}(1)=\left\{
\begin{array}{ll} p \qquad \qquad & \text{if} \,\, k=p^n\,\, \text{for some prime number}\, \, p \\
1\qquad \qquad & \text{if} \,\, k\neq p^n\,\, \text{for any prime
number}\,\, p. \end{array} \right.
\]
\end{lem}
\proof By induction on $k>1$, Lemma  \ref{LL} follows from the
equalities
\[  \Phi
_k(t)=\frac{\displaystyle \sum_{i=0}^{k-1}t^i}{\displaystyle
\prod_{\scriptstyle d\mid k \atop\scriptstyle d<k}\Phi_d(t)},
\quad \Phi _k(1)=\frac{k}{\displaystyle \prod_{\scriptstyle d\mid
k \atop\scriptstyle d<k}\Phi_d(1)},\quad \text{and}\quad \Phi
_k(0)=\frac{1}{\displaystyle \prod_{\scriptstyle d\mid k
\atop\scriptstyle d<k}\Phi_d(0)}.\]  \edvo

Corollary \ref{cor} is a consequence of
\begin{cor} \label{Cor} Let $G$ be an irreducible Hurwitz $C$-group
of degree $m=p^n$, where $p$ is a prime number. Then
\begin{itemize}
\item[($i$)] $\Delta (t)\equiv 1$; \item[($ii$)] the group
$G^{\prime}/G^{\prime \prime}$ is a finite abelian group.
\end{itemize}
\end{cor}

\proof By Theorem \ref{Main}, all roots of $\Delta (t)$ are $m$-th
roots of unity. Let $\lambda$ be one of the roots. Assume that
$\lambda $ is a primitive $p^k$-th root of unity, $1\leq k\leq n$.
Then $\lambda$ is a root of the $p^k$-th cyclotomic polynomial
$$ \Phi _{p^k}(t)= \sum_{i=0}^{p-1}t^{ip^{k-1}} $$
and there is a polynomial $f(t)\in \mathbb Z [t]$ such that
$\Delta (t)=\Phi _{p^k}(t)f(t)$. By Theorem \ref{roots},
$\Phi_{p^k}(1)f(1)=\pm 1$. On the other hand, $f(1)\in \mathbb Z$
and $\Phi_{p^k}(1)=p$. Therefore $\Delta (t)$ has no roots. Thus,
$\deg \Delta (t)=0$ and the group $G^{\prime}/G^{\prime \prime}$
has no free part, that is, it is a finite abelian group. \edvo

\begin{lem} \label{lemma1} Let $G_1$ and $G_2$ be $C$-groups and
$\Delta_1(t)$ and $\Delta_2(t)$ their Alexander polynomials.
Assume that there is a $C$-epimorphism $f:G_1\to G_2$. Then
$\Delta_2(t)$ is a divisor of $\Delta_1(t)$.
\end{lem}

\proof Denote by $N_i$ the kernel of the canonical $C$-epimorphism
$\nu_i:G_i\to \mathbb F_1$, $i=1,2$. It is easy to see that the
homomorphism $g$ in the commutative diagram
\\

\begin{picture}(0,0)(0,0)
\put(100,-5){$1\longrightarrow  N_1  \longrightarrow  G_1
\buildrel{\nu_{1}}\over\longrightarrow  \mathbb F_1
\longrightarrow 1$} \put(100,-50){$1\longrightarrow  N_2
\longrightarrow G_2 \buildrel{\nu_{2}}\over\longrightarrow \mathbb
F_1 \longrightarrow 1$} \put(136,-10){\vector(0,-1){25}}
\put(127,-23){$g$} \put(178,-10){\vector(0,-1){25}}
\put(168,-23){$f$} \put(215,-10){\vector(0,-1){25}}
\put(218,-23){$\simeq$}

\end{picture}
\\
\\
\\
\\
\\
is an epimorphism. This diagram induces the following commutative
diagram
\\

\begin{picture}(0,0)(0,0)
\put(70,-5){$1\longrightarrow  N_1/N'_1  \longrightarrow G_1/N'_1
\buildrel{\nu_{1*}}\over\longrightarrow  \mathbb F_1
\longrightarrow 1$} \put(70,-50){$1\longrightarrow  N_2/N'_2
\longrightarrow G_2/N'_2 \buildrel{\nu_{2*}}\over\longrightarrow
\mathbb F_1 \longrightarrow 1$} \put(121,-10){\vector(0,-1){25}}
\put(111,-23){$g_*$} \put(176,-10){\vector(0,-1){25}}
\put(165,-23){$f_*$} \put(222,-10){\vector(0,-1){25}}
\put(225,-23){$\simeq$} \put(320,-23){$(*)$}

\end{picture}
\\
\\
\\
\\
\\
in which $g_*$ is also an epimorphism.

It follows from diagram $(*)$ that $\Delta_2(t)$ is a divisor of
$\Delta_1(t)$, since $h_2(g_*(n))=g_*(h_1(n))$ for any $n\in
N_1/N'_1$. \edvo

\begin{thm} \label{div} Let $G$ be a Hurwitz $C$-group of degree
$m$. Then its Alexander polynomial $\Delta (t)$ divides the polynomial $(t-1)(t^m-1)^{m-2}$.
\end{thm}
\proof Consider the Hurwitz $C$-group
\begin{equation} \label{uni}
\widetilde G_m=< x_1,\dots,x_m\, \mid \, [x_i,x_1\dots x_m]=1,\,
i=1,\dots,m>.\end{equation}

For any Hurwitz $C$-group $G$ of degree $m$ there is a natural
$C$-epimorphism $f:\widetilde G_m\to G$ sending the $C$-generators
$x_i$ of $\widetilde G_m$ to the $C$-generators $x_i$ of $G$ for
which the product $x_1\dots x_m$ belongs to the center of $G$.
Therefore to prove Theorem \ref{div}, it is sufficient to show
that the Alexander polynomial of $\widetilde G_m$ is equal to
$(-1)^{m-1}(t-1)(t^m-1)^{m-2}$.

Denote by $\widetilde N_m$ the kernel of $\nu : \widetilde G_m\to
\mathbb F_1$ and put  $y=x_1\dots x_m$. Without loss of
generality, we can assume that $\widetilde h(n)=x_mnx_m^{-1}$ for
$n\in \widetilde N_m$. In \cite{Ku.O}, applying the Reidemeister --
Schreier method, it was shown that $\widetilde N_m$ is generated
be the elements
\begin{equation} \label{gen}
a_{k,j}=x_m^kx_jx_m^{-(k+1)}, \end{equation} where $j=2,\dots
,m-1$, $k\in \mathbb Z$, and by the elements
\begin{equation}\label{gen2} a_{k,m} = x_m^kyx_m^{-(k+m)},\end{equation}
where $k\in \mathbb Z$. Then the action $\widetilde h$ is given by
$\widetilde h(a_{k,j})=a_{k+1,j}$.

The relations
$$ yx_j = x_jy, \quad j=2,\dots ,m,$$
give rise (see \cite{Ku.O}) to relations
 \begin{equation} \label{eq2}
a_{k,m} =a_{k+1,m} \end{equation} for $k\in \mathbb Z$ and to relations
\begin{equation}\label{eq1}
a_{k,m}a_{k+m,j}a_{k+1,m}^{-1}=a_{k,j}
\end{equation}
for $j=2,\dots ,m-1$ and $k\in \mathbb Z$. Therefore $\widetilde
N_m$ is a free group generated by $a_{0,m}$ and $a_{k,j}$,
$k=1,\dots, m$, $j=2,\dots, m-1$.

We have $\widetilde h(a_{0,m})=a_{0,m}$ and
$$ \widetilde h(a_{k,j})=a_{k+1,j}$$
for $k=1,\dots, m-1$, $j=2,\dots, m-1$ and
$$\widetilde h(a_{m,j})=a_{0,m}^{-1}a_{1,j}a_{0,m}$$
for $j=2,\dots, m-1$.

Let $\overline a_{k,j}$ be the image of $a_{k,j}$ in $\widetilde
N_m/\widetilde N'_m$. Then the action $h$ is given by
\begin{align*} \widetilde
h(\overline a_{0,m}) & =\overline a_{0,m},  \\
 h(\overline a_{k,j}) & =\overline a_{k+1,j} \quad \text{for} \, \, k=1,\dots,
 m-1, \, j=2,\dots, m-1, \\
 h(\overline a_{m,j}) & =\overline a_{1,j} \qquad \text{for} \, j=2,\dots, m-1 .
 \end{align*}
Simple computations show that the characteristic polynomial of $h$
is equal to $(-1)^{m-1}(t-1)(t^m-1)^{m-2}$. \edvo\\

It is easy to check (\cite{Ku3}, Lemma 4) that $G/G'$ is a
finitely generated free abelian group for any $C$-group $G$.
Moreover, the canonical epimorphism $\text{ab}: G\to G/G'$ is a
$C$-homomorphism if we choose $\text{ab}(x_i)$ as the
$C$-generators of $G/G'$, where the set $\{ x_i\}$ is the set of
$C$-generators of the group $G$. We say that a $C$-group $G$ {\it
consists of $n$ irreducible components} if $G/G'\simeq \mathbb
Z^n$. The notion of the number of irreducible components of a
Hurwitz $C$-group is explained by the following simple lemma.
\begin{lem} \label{number} A (topological) Hurwitz curve $\bar H$
consists of $n$ irreducible components iff its fundamental group
$\pi_1=\pi_1(\mathbb C^2\setminus H)$ consists of $n$ irreducible
components.
\end{lem}

Trivial computations show that the Alexander polynomial $\Delta
(t)$ of an abelian $C$-group $\mathbb Z^n$ is equal to
$(-1)^{n-1}(t-1)^{n-1}$. Therefore, Lemma \ref{lemma1} implies the
following
\begin{lem} \label{le} The Alexander polynomial $\Delta
(t)$ of a Hurwitz $C$-group $G$ consisting of $n$ irreducible
components is divisible by $(t-1)^{n-1}$.
\end{lem}

\begin{thm} \label{mult}
Let a Hurwitz $C$-group $G$ consist of $n$ irreducible components
and $\Delta (t)$ be its Alexander polynomial. Then $$\Delta
(t)=(t-1)^{n-1}P(t),$$ where the polynomial $P(t)\in \mathbb Z[t]$
satisfies $P(1)\neq 0$.
\end{thm}

\proof Let $m$ be the degree of the Hurwitz $C$-group $G$. To obtain a
$C$-presentation of $G$, it is sufficient to add several
$C$-relations to presentation (\ref{uni}).

Since $G$ consists of $n$ irreducible components, the set
$\{1,\dots,m\}$ splits into a disjoint union of $n$ subsets
$J_1,\dots,J_n$ such that $j_1,j_2\in J_k$ iff $x_{j_1}$ and
$x_{j_2}$ are conjugated in $G$. Let $J_k=\{ j_{1,k}<\dots
<j_{m_k,k}\}$. Without loss of generality, one can assume that
among the added $C$-relations there are relations
\begin{equation} \label{added}
x_{j_{i,k}}w_{j_{i,k},j_{i+1,k},1}=w_{j_{i,k},j_{i+1,k},1}x_{j_{i+1,k}}\end{equation}
for $i=1,\dots ,m_k-1$, $k=1,\dots,n$, and for some words
$w_{j_{i,k},j_{i+1,k},1}$. Let
$\nu(w_{j_{i,k},j_{i+1,k},1})=x^{t_{j_{i,k}}}$, where $x$ is the
$C$-generator of $\mathbb F_1$.

Consider diagram (*) in which $G_1=\widetilde G_m$ and $G_2=G$. In
the notations of the proof of Theorem \ref{div}, the elements
$\overline a_{0,m}$ and $\sum_{k=1}^m\overline a_{k,j}$,
$j=2,\dots,m-1$, generate the eigenspace $(N_1/N'_1)_1$ of
$N_1/N'_1$ corresponding to the eigenvalue $1$. Its image
$g_*((N_1/N'_1)_1)\otimes \mathbb C$ is the eigenspace of
$N_2/N'_2\otimes \mathbb C$ corresponding to the eigenvalue
$1$.

Applying the Reidemeister -- Schreier method, the element
$a_{0,m}$, defined in (\ref{gen2}), and the elements $a_{k,j}$,
$k=1,\dots, m$, $j=2,\dots, m-1$, defined in (\ref{gen}) (more
precisely, their images under $g_*$), also generate the group
$N_2$ and each relation (\ref{added}) (after the substitution
$x_1=y(x_2\dots x_m)^{-1}$) gives rise to the relations
\begin{equation} \label{radded1}
a_{r,j_{i,k}}\overline w_{r+1,j_{i,k},j_{i+1,k}}=\overline
w_{r,j_{i,k},j_{i+1,k}}a_{r+t_{j_{i,k}},j_{i+1,k}}\end{equation}
if $1<j_{i,k}<j_{i+1,k}<m$,
\begin{equation} \label{radded2}
a_{0,m}(\prod_{s=2}^{m-1}a_{r+m-s,m+1-s}^{-1}) \overline
w_{r+1,1,j_{i+1,k}}=\overline
w_{r,1,j_{i+1,k}}a_{r+t_{1,j_{i+1,k}}}
\end{equation} if
$1=j_{i,k}<j_{i+1,k}<m$,
\begin{equation} \label{radded3}
a_{r,j_{i,k}}\overline w_{r+1,j_{i,k},m}=\overline
w_{r,j_{i,k},m}\end{equation} if $1<j_{i,k}<j_{i+1,k}=m$, and
\begin{equation} \label{radded4}
a_{0,m}(\prod_{s=2}^{m-1}a_{r+m-s,m+1-s}^{-1}) \overline
w_{r+1,1,m}=\overline w_{r,1,m}\end{equation} if
$1=j_{i,k}<j_{i+1,k}=m$, where $r\in \mathbb Z$ and each word
$\overline
w_{r,j_{i,k},j_{i+1,k}}=x_m^rw_{j_{i,k},j_{i+1,k},1}x_m^{-(r+t_{j_{i,k}})}$
is written in the generators $a_{l,s}$.

As in \cite{Ku.O}, one can show that the words $\overline
w_{r,j_{i,k},j_{i+1,k}}$ and $\overline w_{r+m,j_{i,k},j_{i+1,k}}$
are conjugated in $G_2$ by $a_{0,m}$. Therefore, taking the sum
over $r$, relations (\ref{radded1}) -- (\ref{radded4}) give rise
in $N_2/N'_2$ to the following relations:
\begin{equation} \label{radded11}
\sum_{r=1}^m\overline a_{r,j_{i,k}}=\sum_{r=1}^m\overline
a_{r,j_{i+1,k}}\end{equation} if $1<j_{i,k}<j_{i+1,k}<m$,
\begin{equation} \label{radded21}
m\overline a_{0,m}-\sum_{s=2}^{m-1}\sum_{r=1}^m\overline a_{r,s}
=\sum_{r=1}^m\overline a_{r,j_{i+1,k}}
\end{equation} if
$1=j_{i,k}<j_{i+1,k}<m$,
\begin{equation} \label{radded31}
\sum_{r=1}^m\overline a_{r,j_{i,k}}=0\end{equation} if
$1<j_{i,k}<j_{i+1,k}=m$, and
\begin{equation} \label{radded41}
m\overline a_{0,m}-\sum_{s=2}^{m-1}\sum_{r=1}^m\overline a_{r,s}
=0
\end{equation}
if $1=j_{i,k}<j_{i+1,k}=m$.

It is easy to see that for any decomposition  $\{ 1,\dots, m\}=
\bigsqcup_{k=1}^nJ_k$, equations (\ref{radded11}) --
(\ref{radded41}) are linear independent over $\mathbb Z$, their
number being $m-n$, and the elements $\overline a_{0,m}$ and
$\sum_{r=1}^m\overline a_{r,s}$, $s=2,\dots, m-1$ generate the
eigenspace $(N_1/N'_1)_1$. Therefore the rank of the kernel of the
restriction of $g_{*}$ to $(N_1/N'_1)_1$ is not less then $m-n$.
Thus $$\dim (N_2/N'_2\otimes \mathbb C)_1\leq \dim
(N_1/N'_1\otimes \mathbb C)_1-(m-n)
=n-1.$$
On the
other hand, by Lemma \ref{le}, the Alexander polynomial $\Delta
(t)$ of a $C$-group $G$ consisting of $n$ irreducible components
is divisible by $(t-1)^{n-1}$. \edvo

\begin{cor} \label{sign}
Let a Hurwitz $C$-group $G$ consist of $n$ irreducible components
and $\Delta (t)$ be its Alexander polynomial. Then
$$\Delta(0)=(-1)^{\deg \Delta (t)-(n-1)} .$$
\end{cor}
\proof The Alexander polynomial
$$\Delta (t)=\det (h_{\mathbb C}
-t\text{Id})=(-t)^{\deg \Delta (t)}+\sum_{i=0}^{\deg \Delta
(t)-1}c_it^i$$
satisfies, and by Theorem \ref{mult}, $\Delta
(t)=(t-1)^{n-1}P(t)$ where the polynomial $P(t)\in \mathbb Z[t]$
is such that $P(1)\neq 0$. Therefore the polynomial
$P(t)=(-1)^{\deg \Delta (t)} \prod_i \Phi_{n_i}(t)$ is a product
(up to sign) of some cyclotomic polynomials $\Phi_{n_i}(t)$ with
$n_i>1$. By Lemma \ref{LL}, $\Phi_{n_i}(0)=1$ for all $n_i>1$.
Therefore $\Delta (0)=(-1)^{n-1}(-1)^{\Delta (t)}=(-1)^{\deg
\Delta (t)-(n-1)}$.  \edvo

\begin{lem} \label{addgenerators}
Let $j:\{ 1,\dots, n\}\to \{ 1,\dots, m\}$ be an injective
function. Assume that a $C$-group $G$ is generated by
$C$-generators $x_1,\dots,x_n$ and $w=\overline x_1\dots \overline
x_m$ is a quasipositive word in $x_1,\dots,x_n$ {\rm (}that is,
each $\overline x_k$ is conjugated to some $x_{i_k}\in
\{x_1,\dots,x_n\}${\rm )} such that $\overline x_{j(i)}=x_i$ for
$i=1,\dots, n$. If $w$ belongs to the center of $G$, then $G$ is a
Hurwitz $C$-group of degree $m$.
\end{lem}

\proof Let 
$G=<x_{1},\dots,x_{n}\, \mid \, R>$ 
be a $C$-presentation of the group $G$. Put $J=\{ 1\leq j\leq m \,
\mid \, j=j(i),\, i=1,\dots, n\}$. By assumption, we have
$\overline x_{j(i)}=x_i$ and
\begin{equation} \label{pre2} \overline x_j
=w_j^{-1}x_{i_j}w_j\end{equation} for $j\not\in J$, where $w_j$ is
a word in $x_1,\dots,x_n$. Note that relations (\ref{pre2}) are
$C$-relations. Therefore if we add the generators $\overline x_j$,
$j\not\in J$, to the set of generators $\{ x_1 ,\dots,x_n\}$ and
add relations (\ref{pre2}) to $R$, then we obtain a
$C$-presentation of the same group $G$. To complete the proof, it
suffices to renumber the obtained set of generators. \edvo

\begin{prop} \label{prod}
Let $G_i$, $i=1,2$, be a Hurwitz $C$-group of degree $m_i$ and
$\Delta_i(t)$ its Alexander polynomial. Then there exists a
Hurwitz $C$-group of degree $2m_1m_2$ with Alexander
polynomial $\Delta(t)=\Delta_1(t)\Delta_2(t)$.
\end{prop}

\proof Let $G_i=<x_{1,i},\dots,x_{m_i,i}\, \mid \, R_i > $ be a
Hurwitz $C$-presentation of the group $G_i$ of degree $m_i$.

Consider the amalgamated product $\widetilde G=G_1*_{\{
x_{m_1,1}=x_{m_2,2}\}}G_2$. It is a $C$-group  given by the
presentation
$$\widetilde G=<x_{1,i},\dots,x_{m_i,i},\, i=1,2\, \mid \, R_1\cup R_2, \, x_{m_1,1}=x_{m_2,2} > .$$
Put $y_i=x_{1,i}\dots x_{m_i,i}$, $i=1,2$, and denote by $N_i$
(resp. $\widetilde N$) the kernel of $\nu: G_i\to \mathbb F_1$
(resp. of $\nu: \widetilde G\to \mathbb F_1$).

As in the proof of Theorem \ref{div}, applying the Reidemeister --
Schreier method, one can show that the group $N_i$ (resp.
$\widetilde N$) is generated by the elements
$a_{k,j,i}=x_{m_i,i}^kx_{j,i}x_{m_i,i}^{-(k+1)}$, where $j=2,\dots
,m_i-1$, $k\in \mathbb Z$, and by the elements $ a_{k,m_i,i} =
x_{m_i,i}^ky_ix_{m_i,i}^{-(k+m_i)}$, $k\in \mathbb Z$ (resp. by
the union of these elements, since $x_{m_1,1}=x_{m_2,2}$ in
$\widetilde G$). The set of defining relations of $N_i$ (resp. of
$\widetilde N$) is obtained from the set $\overline R_i=\{
x_{m_i,i}^nr_{l,i}x_{m_i,i}^{-n} \mid r_{l,i}\in R_i\, n\in
\mathbb Z\}$ (resp. from $\overline R=\overline R_1\cup\overline
R_2$) after rewriting the words in the alphabet $\{ a_{k,j,i}\}$
(resp. in the union of these alphabets). Therefore $\widetilde
N=N_1*N_2$ is the free product of the groups $N_1$ and $N_2$.

Note that it follows from the proof of Theorem \ref{div} that the
elements $a_{k,j,i}$ and $a_{k+lm_i,j,i}$ are conjugated in $N_i$
for all $l\in\mathbb Z$.

Let $\widetilde h_i$ be the automorphism of $N_i$ given by
conjugation by $x_{m_i,i}$. Then the automorphism $\widetilde h$
of $\widetilde N$ given by conjugation by $x_{m_1,1}=x_{m_2,2}$,
which is equal to $\widetilde h_1*\widetilde h_2$. Therefore the
Alexander polynomial $\widetilde \Delta (t)$ of the group
$\widetilde G$ is equal to $\widetilde
\Delta(t)=\Delta_1(t)\Delta_2(t)$.

Consider a group
\begin{align*}  G= & <x_{1,i}\dots ,x_{m_i,i}, i=1,2 \mid R_1\cup
R_2,\, x_{m_1,1}=x_{m_2,2}, \\ & [x_{j,i},y_{\overline i}^{m_i}]=1
, j=1,\dots,m_i-1, i=1,2> ,\end{align*} where $\overline i=\{
1,2\}\setminus \{ i\}$ (recall that $x_{1,\overline i},\dots,
x_{m_{\overline i},\overline i}$ commute with $y_{\overline i}$).
Let $N$ be the kernel of $\nu :G\to \mathbb F_1$.

To obtain a presentation of the group $N$ from the presentation of
the group $\widetilde N$ described above, one should add the
relations induced by the relations
$$[x_{j,i},y_{\overline i}^{m_i}]=1
, j=1,\dots,m_i-1, i=1,2.$$ It is easy to see that these
additional relations are
\begin{equation} \label{N}
a_{k,j,i}a_{0,m_{\overline i}}^{m_i}=a_{0,m_{\overline
i}}^{m_i}a_{k+m_1m_2,j,i},\end{equation} since $a_{k,m_{\overline
i}}=a_{k+1,m_{\overline i}}$ in $N_{\overline i}$ for all $k$,
that is, the additional relations (\ref{N}) say that $a_{k,j,i}$ and
$a_{k+m_1m_2,j,i}$ are conjugated. But these elements were
conjugated in $\widetilde N$. Therefore $\widetilde N/\widetilde
N'\simeq N/N'$ and the groups $\widetilde G$ and $G$ have the same
Alexander polynomials. To complete the proof of Proposition
\ref{prod}, let us notice that $y_1^{m_2}y_2^{m_1}=(x_{1,1}\dots
x_{m_1,1})^{m_2}(x_{1,2}\dots x_{m_2,2})^{m_1}$ belongs to the
center of the group $G$. Therefore, by Lemma \ref{addgenerators},
$G$ is a Hurwitz $C$-group of degree $2m_1m_2$.  \edvo\\

For two Hurwitz $C$-groups $G_1$ and $G_2$ given by Hurwitz
$C$-presenta\-tions $G_i=<x_{1,i},\dots,x_{m_i,i}\, \mid \, R_i
> $, $i=1,2$,  the Hurwitz $C$-group $G$, constructed in the proof of
Proposition \ref{prod}, is called {\it a Hurwitz product} of $G_1$
and $G_2$. A Hurwitz product of $G_1$ and $G_2$ will be denoted by
$G_1\diamond G_2$. Of course, a Hurwitz product of $G_1$ and $G_2$
depends on Hurwitz $C$-presentations of $G_1$ and $G_2$, but by
Proposition \ref{prod}, the Alexander polynomial of $G_1\diamond
G_2$ does not depend on the Hurwitz $C$-presentations of the
factors.

\begin{lem} \label{linkgroups} The fundamental group $G_{n,m}=\pi_1(\mathbb C^2\setminus
C_{n,m})$ of the complement of the affine plane algebraic curve
$C_{m,n}$, given by the equation $w^n-z^m=0$, where $n$ and $m$
are any positive integers, is a Hurwitz $C$-group.
\end{lem}
{\bf Remark}. Note that this lemma does not follow from the
statement, mentioned above, on the fundamental group of the
complement of an affine Hurwitz curve, since it is assumed there
that the line at infinity is in general position with respect to
the Hurwitz curve. Here, the line at infinity is in special position. If we consider the local fundamental group
$G=\pi_1(B_{\varepsilon}\setminus C)$, where $C$ is an irreducible
singularity in a small ball $B_{\varepsilon}$, then $G$ has always
a natural structure of an irreducible $C$-group. It has a
non-trivial center iff the singularity $C$ is of type $x^p=y^q$
with $p$ and $q$ coprime (see \cite{Bu-Zi}).

 \proof Indeed, a braid
monodromy of the singularity $w^n=z^m$ with respect to the
projection $(z,w)\to z$ is equal to
$$b_{n,m}=(\sigma_1\dots \sigma_{n-1})^m\in \text{Br}_n,$$
where $\sigma_1,\dots ,\sigma_{n-1}$ are standard generators of
the braid group $\text{Br}_n$, that is, the generators satisfy the
following relations
$\sigma_i\sigma_{i+1}\sigma_i=\sigma_{i+1}\sigma_{i}\sigma_{i+1}$
for $i=1,\dots,n-2$ and $[\sigma_i,\sigma_j]=1$ for $\mid i-j\mid
\geq 2$. The group $\text{Br}_n$ acts on the free group $\mathbb
F_n$ generated by $x_1,\dots,x_n$. This action is given by
$\sigma_j(x_i)=x_i$ if $j\neq i,i+1$, $\sigma_j(x_{j+1})=x_j$, and
$\sigma_j(x_{j})=x_jx_{j+1}x_j^{-1}$. Denote by $B_{n,m}$ the
cyclic subgroup of $\text{Br}_n$ generated by $b_{n,m}$. Then (see
\cite{Ku}) the group
$$G_{n,m}=<x_1,\dots,x_n\,\mid \, x_i=b(x_i), i=1,\dots,n,\, \,
b\in B_{n,m}>$$ is a $C$-group, and by Lemma \ref{addgenerators},
it is a Hurwitz $C$-group, since $b_{n,m}^n=(\Delta^2_n)^m\in
B_{n,m}$, where $\Delta_n$ is the Garside element of the braid
group $\text{Br}_n$, and therefore the element $(x_1\dots x_n)^m$
belongs to the center of $G_{n,m}$ (since
$\Delta^2_n(x_i)=(x_1\dots x_n)x_i(x_1\dots x_n)^{-1}$ and
$\Delta^2_n(x_1\dots x_n)=x_1\dots x_n$). \edvo
\begin{prop}
\label{Le} {\rm (\cite{L}) } If $m$ and $n$ are coprime, then the
group $G_{n,m}$ has the Alexander polynomial 
$$ \Delta_{n,m}(t)=\frac{(t-1)(t^{nm}-1)}{(t^n-1)(t^m-1)}.$$
\end{prop}
\begin{prop} \label{2m} The Alexander polynomial of the group $G_{2,2m}$ is equal to
$$\Delta(t)=(1-t)\sum_{i=0}^{m-1}t^{2i}.$$
\end{prop}
\proof It follows from the proof of Lemma \ref{linkgroups} that

$$G_{2,2m}=<x_1,x_{2}\, \mid \,
[x_1,(x_1x_2)^m]=[x_2,(x_1x_2)^m]=1>. $$

Let us show that the relations
\begin{equation} \label{commut}
[x_2,(x_1x_2)^m]=1, \, \, \, \, i=1,2, \end{equation} are
equivalent to the single relation $$(x_1x_2)^m=(x_2x_1)^m.$$
Indeed, relations (\ref{commut}) imply

$$(x_1x_2)^m=x_1^{-1}(x_1x_2)^mx_1=(x_2x_1)^m,$$

the relation $(x_1x_2)^m=(x_2x_1)^m$ implies

$$ x_2(x_1x_2)^m=(x_2x_1)^mx_2=(x_1x_2)^mx_2$$ and

$$ x_1(x_1x_2)^m=x_1(x_2x_1)^m=(x_1x_2)^mx_1.$$
Therefore
$$G_{2,2m}=<x_1,x_{2}\, \mid \,
(x_1x_2)^m(x_2x_1)^{-m}=1>. $$ Denote by
$r=(x_1x_2)^m(x_2x_1)^{-m}$. Applying the free differential
calculus of Fox (\cite{C-F}), it is easy to see that
$$\nu_*\bigl(\frac{\partial r}{\partial x_1}\bigr)=
1+t^2+\dots + t^{2(m-1)}-t^{2m-1}-t^{2m-3}-\dots -t^1$$ and
$$\nu_*\bigl(\frac{\partial
r}{\partial x_2}\bigr)= t+\dots +
t^{2m-1}-t^{2(m-1)}-t^{2(m-2)}-\dots -t^2-1.$$ Therefore
$$\Delta(t)=(1-t)\sum_{i=0}^{m-1}t^{2i}.$$ \edvo\\

Consider the group \begin{equation} \label{G2} \begin{array}{rl}
G(2)=<x_1,\dots,x_4\mid  &
x_2^2x_1x^{-2}_2=x_4,\, x_3=x_2, \, \, x_4^2x_2x_4^{-2}=x_2 \\
 & [x_i,x_1\dots
x_4]=1\,\,\text{for}\, \, i=1,\dots, 4>.\end{array}
\end{equation}
\begin{prop} \label{2}
The Alexander polynomial of the group $G(2)$ is equal to
$\Delta(t)=t^2-1.$
\end{prop}
\proof Denote by $N(2)$ the kernel of $\nu: G(2)\to \mathbb F_1$
and put $m=4$ and $y=x_1\dots x_4$. In the notations of the proof
of Theorem \ref{div}, it follows from the relations $[x_i,y]=1$
for $i=1,\dots,4$ that the group  $N(2)$ is generated by the
elements $a_{k,j}=x^k_{4}x_jx_{4}^{-(k+1)}$, $k=1,\dots, 4$,
$j=2,3$, and the element $a_{0,4}=yx_{4}^{-4}$. In our case
relations
(\ref{eq2}) and (\ref{eq1}) have the form 
\begin{equation} \label{0de}
a_{k,4}=a_{0,4}\qquad \text{and}\qquad
a_{k+4,j}=a_{0,4}^{-1}a_{k,j}a_{0,4}\end{equation} for all $k$.

The relation $x_3=x_{2}$ gives rise to the relations
\begin{equation} \label{1de}
a_{k,3}=a_{k,2}
\end{equation}
for all $k$.

The relation $x_{4}^2x_2=x_2x_{4}^2$ gives rise to the relations
\begin{equation} \label{2de}
a_{k+2,2}=a_{k,2}
\end{equation}
for all $k$.

 The relation $x^2_2x_{1}=x_{4}x^2_2$,
written as $x^2_2y=x_{4}x_2^4x_4$ (since $x_2=x_3$), gives rise to
the relations
\begin{equation} \label{4de}
a_{k,2}a_{k+1,2}a_{k+2,4}=a_{k+1,2}a_{k+2,2}a_{k+3,2}a_{k+4,2}
\end{equation}
for all $k$.

It follows from (\ref{0de}) -- (\ref{4de}) that $N(2)$ is
generated by $a_{1,2}$, $a_{2,2}$ and $a_{0,4}$, being subject to
the relations
$$a_{0,4}=(a_{1,2}a_{2,2})^{-1}(a_{2,2}a_{1,2})^2=(a_{2,2}a_{1,2})^{-1}(a_{1,2}a_{2,2})^2$$
and
$$ [a_{1,2},a_{0,4}]=[a_{2,2},a_{0,4}]=1.$$
Therefore the group $N(2)/N(2)'$ is a free abelian group generated
by the images $\overline a_{1,2}$ and $\overline a_{2,2}$ of the
elements $a_{1,2}$ and $a_{2,2}$.

As in the proof of Theorem \ref{div}, the action $\widetilde h$ on
$N(2)$ is given by $\widetilde h(a_{1,2})=a_{2,2}$, $\widetilde
h(a_{2,2})=a_{3,2}=a_{1,2}$. The induced action $h$ on
$N(2)/N(2)'$ is given by $h(\overline a_{1,2})= \overline a_{2,2}$
and  $h(\overline a_{2,2})= \overline a_{1,2}$  the characteristic
polynomial of which is equal to $(t-1)(t+1)$. \edvo
\begin{cor} For any $k\in \mathbb N$ there is a
Hurwitz $C$-group $G$ consisting of two irreducible components whose Alexander polynomial $\Delta (t)=(t-1)P(t)$ is such that
$\mid P(1)\mid =k$. \end{cor}

\proof If $k>2$ then, by Proposition \ref{2m}, we have $P(1)=-k$,
where $P(t)=-\sum_{i=0}^{k-1}t^{2i}$ is the factor of the
Alexander polynomial $\Delta(t)=(1-t)\sum_{i=0}^{k-1}t^{2i}$ of
the group $G_{2,2k}$. If $k=2$, then, by Proposition \ref{2}, the
group $G(2)$ has the desired property, since its Alexander
polynomial is $\Delta (t)=(t-1)(t+1)$. In the case $k=1$, one can
take the abelian Hurwitz $C$-group $G=\mathbb Z^2$. \edvo

\begin{prop} \label{irr}
For any $k\in \mathbb N$ there exists
\begin{itemize}
\item[($i$)] an irreducible Hurwitz $C$-group whose Alexander
polynomial $\Delta (t)$ has $\deg \Delta (t)=2k$; \item[($ii$)] a
Hurwitz $C$-group consisting of two irreducible components and
whose Alexander polynomial $\Delta (t)=(t-1)P(t)$ satisfies
$\deg P(t)=k$.
\end{itemize}
\end{prop}
\proof By Propositions \ref{prod} and \ref{Le}, the Alexander
polynomial $\Delta (t)$ of a Hurwitz product $G_{2,3}^{\diamond
k}$ is equal to $(t^2-t+1)^k$.

To prove ($ii$), it suffices to take the groups $G(2)\diamond
G_{2,3}^{\diamond n}$ if $k=2n+1$ is odd and $\mathbb Z^2\diamond
G_{2,3}^{\diamond n}$ if $k=2n$ is even. \edvo

\begin{que} Let $P(t)\in \mathbb Z[t]$ be a polynomial whose roots are roots
  of unity, let
$t=1$ be a root of $P(t)$ of multiplicity $k$, and
$P(0)=(-1)^{\deg P(t)-k}$. Assume also that $P(1)=1$ if $k=0$. Does there exist a Hurwitz $C$-group $G$ with Alexander
polynomial $\Delta (t)=P(t)${\rm ?}
\end{que}

\section{The first Betti number of cyclic coverings of the plane}
\label{section6}

Consider the infinite cyclic covering $f =f_{\infty }:X_{\infty }
\to X'={\Bbb C}^{2}\, \setminus \, H$ corresponding to the
epimorphism $\nu :\pi_1(\mathbb C^2\setminus H)\to \mathbb F_1$.
Let $h\in \text{Deck}(X_{\infty } /X')\simeq {\Bbb F}_{1}$ be a
covering transformation corresponding to the $C$-generator $x\in
{\Bbb F}_{1}$. We say that $h$ is the {\it monodromy of a Hurwitz
curve} $H$. The space $X'$ will be considered as the quotient
space $X'=X_{\infty }/{\Bbb F}_{1}$. In such a situation Milnor
\cite{Mi} considered an exact sequence of chain complexes
$$ 0\to C_{\cdot }(X_{\infty })\buildrel{h-id}\over\longrightarrow
 C_{\cdot }(X_{\infty })\buildrel{\varphi _{*}}\over\longrightarrow
 C_{\cdot }(X')\to  0 $$
which gives an exact homology sequence \begin{equation}
\label{chain} \ldots \to H_{1}(X_{\infty
})\stackrel{h-id}\longrightarrow H_{1}(X_{\infty })\stackrel{f
_*}\longrightarrow H_{1}(X')\rightarrow H_{0}(X_{\infty
})\rightarrow  0 \end{equation} (We often write $h$ instead of
$h_{*}$, if it does not lead to a misunderstanding).

If $G_{n}\subset {\Bbb F}_{1}$ is an infinite cyclic group
generated by $h^{n}$, then $X'_n= X_{\infty }/G_{n}$ and $X'=
X'_n/\mu _{n}$, where $\mu _{n}={\Bbb F}_{1}/G_{n}$ is the cyclic
group of order $n$. Denote by $h_{n}$ an automorphism of $X'_n$
induced by the monodromy $h$. Then $h_{n}$ is a generator of the
covering transformation group $\text{Deck}(X'_{n} /X')=\mu _{n}$
acting on $X'_n$. We apply the sequence
\begin{equation}
\label{chain-n} \ldots \to H_{1}(X_{\infty
})\stackrel{h^n-id}\longrightarrow H_{1}(X_{\infty
})\stackrel{g_{n,*}}\longrightarrow H_{1}(X'_n)\rightarrow
H_{0}(X_{\infty })\rightarrow  0 \end{equation}
 constructed in the same way as (\ref{chain})
to the infinite cyclic covering $g_n=g_{\infty ,n}:X_{\infty }\to
X'_n$, to analyse the group $H_{1}(X'_{n})$.

Denote by $H_1(X_{\infty}, \mathbb C)_n$ the subspace of
$H_1(X_{\infty},\mathbb C)$ corresponding to the eigenvalues
$\lambda$ of $h_*$ which are $n$-th roots of unity and denote by
$H_1(X_{\infty}, \mathbb C)_{n,\neq 1}$ the subspace corresponding
to the eigenvalues $\lambda\neq 1$ of $h$ which are $n$-th roots
of unity. Obviously, $\dim H_1(X_{\infty}, \mathbb C)_n= r_n$ and
$\dim H_1(X_{\infty}, \mathbb C)_{n,\neq 1}=r_{n,\neq 1}$, where
$r_n$ (resp. $r_{n,\neq 1}$) is the number of roots of the
Alexander polynomial $\Delta (t)$ of the Hurwitz curve $\bar H$
which are $n-th$ roots of unity (resp. not equal to $1$). Note
that by Lemma \ref{number} and Theorem \ref{mult},
$$r_n-r_{n,\neq
1}=r_1=\# \{ \text{irreducible components of} \, \bar H\} -1.$$

\begin{prop} \label{dim-X-n}
We have \begin{itemize} \item[($i$)] $\dim H_1(X'_{n}, \mathbb C)=
r_n+1,$ \item[($ii$)] $\dim H_1(X'_{n}, \mathbb C)_1= r_1+1=\# \{
\text{irreducible components of} \, \bar H\}.$
\end{itemize}
\end{prop}

\proof This follows from the exact sequence (\ref{chain-n}). \edvo

Let $\bar{H}$ be a Hurwitz curve consisting  of $k$ irreducible components
$\bar H_1,\dots, \bar H_k$. Choose a line $L\subset {\Bbb C}^{2}$
belonging to the pencil of lines (with respect to which $\bar H$
is defined) and transversely intersecting the curve $H$. Denote
by $\gamma _{i}$ a circle of small radius in $L$ with center
at one of the intersection points $H_{i}\cap L$. It is easy to see
that the cycles $\gamma _{1},\, \dots ,\gamma _{k}$, corresponding
to the chosen loops, form a basis in $H_{1}({\Bbb C}^{2}\,
\setminus \, H,{\Bbb Z})$ and are independent of the choice of the
line $L$. Let $\bar \gamma _{i},\, i = 1,\ldots ,k$, be a cycle in
$H_{1}(X'_n,{\Bbb Z})$ corresponding to a simple path $f ^{-1}_{n}
(\gamma_i)$.
\begin{lem} The cycles $\bar \gamma _{i},\, i = 1,\ldots ,k$, are linear
independent in $H_{1}(X'_n,{\Bbb Z})$ and form a basis in
$H_1(X'_n)_1$.
\end{lem}

\proof Obviously, all $\bar \gamma_i$ are invariant under the
action $h_n$. Now the proof follows from Proposition \ref{dim-X-n}
($ii$) and from the remark that under the homomorphism
$(f_{n})_{*}:H_{1}(X'_n,{\Bbb Z})\to H_{1}({\Bbb C}^{2}\,
\setminus \, H,{\Bbb Z})$ we have $(f_{n})_{*}(\bar \gamma
_{i})=n\gamma _{i}$. \edvo\\

In the notations of the proof of Theorem \ref{sympl-thm}, the
covering $f_n$ can be extended to a map ${\widetilde
f}_n:{\widetilde X}_n\to \mathbb C\mathbb P^2$ branched along
$\bar H$ and, maybe, along $\overline L_{\infty}$. Here
$\widetilde X_n$ is a closed four dimensional variety locally
isomorphic over a singular point of $\bar H$ to a complex analytic
singularity given by an equation $w_1^n=\widetilde F_1(u_1,v_1)$,
where $\widetilde F_1(u_1,v_1)=\prod(v_1-v_{1,j}(u_1))$ and the
product is taken over those branches of $\bar H$ for
which the closure contains the singular point of $\bar H$. In addition,
$\widetilde X_n$ is locally isomorphic over a neighbourhood of an
intersection point of $\bar H$ and $L_{\infty}$ to the singularity
given by $w_2^n=(v_2-e^{\frac{2\pi i}{m}})u^{d}$, where $d$ is the
smallest non-negative integer for which $m+d$ is divisible by $n$.
The variety $\widetilde X_n$, if $\widetilde
f_n^{-1}(L_{\infty})\subset \text{Sing}\, \widetilde X_n$, can be
normalized (as in the algebraic case) and we obtain a covering
$\widetilde f_{n,\, \text{norm}}:\widetilde X_{n,\,
\text{norm}}\to\mathbb C\mathbb P^2$ in which $\widetilde X_{n,\,
\text{norm}}$ is a singular analytic variety at its finitely many
singular points. One can resolve them and obtain a smooth manifold
$\overline X_n$. Let $\sigma :\overline X_n\to \widetilde X_{n,\,
\text{norm}}$ be a resolution of the singularities,
$E=\sigma^{-1}(\text{Sing}\, \widetilde X_{n,\, \text{norm}})$,
and $\overline f_n=\widetilde f_{n,\, \text{norm}}\circ \sigma$.
Denote by $R_i=\widetilde f_{n,\, \text{norm}}^{-1}(\bar H_i)$,
$i=1,\dots, k$, and $R_{\infty}=\widetilde f_{n,\,
\text{norm}}^{-1}(L_{\infty})$. Note that the restriction of
$\widetilde f_{n,\, \text{norm}}$ to each $R_i$, $i=1,\dots,k$, is
one-to-one and the restriction of $\widetilde f_{n,\,
\text{norm}}$ to $R_{\infty}$ is a $n_0$-sheeted cyclic covering,
where $n_0=\text{GCD}(n,d)$ and the ramification index of
$\widetilde f_{n,\, \text{norm}}$ along $R_{\infty}$ is equal to
$n_{\infty}=\frac{n}{n_0}$. As in the algebraic case, it is easy
to show that $R_{\infty}$ is irreducible. Denote by $\overline
R_i=\sigma^{-1}(R_i)$, $i=1,\dots,k,\infty$, the proper transform
of $R_i$.

We have the embeddings $i_1:X'_n\hookrightarrow X_n=\overline
X_n\setminus E$ and $i_2:X_n\hookrightarrow \overline X_n$.
\begin{lem} \label{ker} The induced homomorphism
$i_{1*}:H_1(X'_n)\to H_1(X_n)$ is an epimorphism with $\ker
i_{1*}=H_1(X'_n)_1$.
\end{lem}
\proof We have $$ X'_n=X_n\setminus (\cup_{i=1}^k\overline
R_i)\cup \overline R_{\infty}$$ and each $\overline R_i$,
$i=1,\dots, k,\infty$, is a codimension two submanifold of $X_n$.
Therefore each 1-dimensional cycle $\gamma\subset X_n$ can be
moved outside of $(\cup_{i=1}^k\overline R_i)\cup \overline
R_{\infty}$. Thus, $i_{1,*}$ is an epimorphism.

Let a complex line $L\subset \mathbb C\mathbb P^2$ meet $L_{\infty}$
transversely  at $q\in L_{\infty}\setminus \bar H$
and $\gamma_{\infty}$ be a simple small loop around $L_{\infty}$
lying in $L$. Then $f_n^{-1}(\gamma_{\infty})$ splits into the
disjoint union of $n_0$ simple loops $\bar \gamma_{\infty,i}$,
$i=1,\dots,n_0$. Since $R_{\infty}$ is irreducible, each two loops
$\bar \gamma_{\infty,i}$ and $\bar \gamma_{\infty,j}$ belong to
the same homology class in $H_1(X'_n)$ (denote it by $\bar
\gamma_{\infty}$). Therefore $n_0\bar\gamma_{\infty}\in
H_1(X'_n)_1$. Now lemma follows from the remark that $\bar
\gamma_1,\dots, \bar\gamma_k, \bar\gamma_{\infty}$ generate $\ker
i_1$ and $\bar \gamma_1,\dots, \bar\gamma_k$ generate
$H_1(X'_n)_1$. \edvo

\begin{lem} \label{isom} The homomorphism $i_{2*} :H_1(X_n,\mathbb C)\to H_1(\overline
X_n,\mathbb C)$ is an isomorphism.
\end{lem}

\proof We have $X_n=\overline X_n\setminus E$. Denote by $T\subset
\overline X_n$ a closed regular neighbourhood of $E$ and let
$\partial T$ be its boundary, $T'=T\setminus E$, and
$T^0=T\setminus \partial T$. It is known (see, for example, the
proof of Proposition 3.4 from \cite{Di}) that the homomorphism
$i_*:H_1(\partial T,\mathbb C)\to H_1(T,\mathbb C)$, induced by
the imbedding $i:\partial T\hookrightarrow T$, is an isomorphism
and, besides, there is a deformation retract $T'\searrow
\partial T$. Therefore there is a deformation retract
$X_n\searrow X_n^0$, where $X_n^0=X_n\setminus T^0$. Now the lemma
follows from the Mayer -- Vietories sequence
$$ H_2(\overline X_n)\to H_1(\partial T)\to H_1(T)\oplus
H_1(X_n^0)\to H_1(\overline X_n)\to 0.$$ \edvo

The proof of Theorem \ref{th} follows from Lemmas \ref{ker},
\ref{isom} and Proposition \ref{dim-X-n}.

\begin{prop} \label{b1}
For any $k\in \mathbb N$, there exists
\begin{itemize}
\item[($i$)] an irreducible Hurwitz curve $\bar H_{k}$ such that a
resolution of singularities $\overline X_{k,6}$ of the  cyclic
covering of $\mathbb C\mathbb P^2$ of degree six, branched along
$\bar H_{k}$, has first Betti number $b_1(\overline
X_{k,6})=2k$; \item[($ii$)] a Hurwitz curve $\bar H_k$ consisting
of two irreducible components such that the first Betti number
$b_1(\overline X_{k,6})$ of a resolution of singularities
$\overline X_{k,6}$ of the cyclic covering of degree six, branched
along $\bar H$, is equal to $k$.
\end{itemize}
\end{prop}
\proof In the proof of Proposition \ref{irr}, it was shown that
the Alexander polynomial $\Delta (t)$ of a Hurwitz product
$G_{2,3}(k)=G_{2,3}^{\diamond k}$ is equal to $(t^2-t+1)^k$ and
that the Alexander polynomials $\Delta (t)$ of  Hurwitz products
$G_{2,3}(2,n)=G(2)\diamond G_{2,3}^{\diamond n}$ and
$G_{2,3}(\text{ab},n)=\mathbb Z^2\diamond G_{2,3}^{\diamond n}$
respectively are equal to $(t-1)(t+1)(t^2-t+1)^n$ and
$(1-t)(t^2-t+1)^n$.

The groups $G_{2,3}(k)$, $G_{2,3}(2,n)$, and
$G_{2,3}(\text{ab},n)$ are Hurwitz $C$-groups. Moreover, one can
assume that the degrees of these Hurwitz $C$-groups are divisible
by six (one can take $y^6$, where $y$ is the product of the
$C$-generators of a Hurwitz $C$-presentation of a group, and apply
Lemma \ref{addgenerators}). Therefore by Theorem 6.2 from
\cite{Ku}, each of these groups can be realized as the fundamental
group $\pi_1(\mathbb C^2\setminus H)$ for some Hurwitz curve of
degree divisible by six. The curve $\bar H$ is irreducible in the
case of $G_{2,3}(k)$ and consists of two irreducible components in
the other two cases. Now Theorem \ref{th} implies Proposition
\ref{b1}. \edvo

\begin{prop} \label{b1G2} For any $k\in \mathbb N$, there is a Hurwitz curve $\bar
H_k$ which consists of $k+1$ irreducible components, has
singularities of the form $w^{2^{3k-1}}-z^{2^{3k-1}}=0$, and which
is the branch curve of a $2$-sheeted cyclic covering $\overline
f_2: \overline X_{k,2}\to \mathbb C\mathbb P^2$ with
$b_1(\overline X_{k,2})=k$.

In particular, the Hurwitz curve $\bar H_1$ has $\deg \bar
H_1=2^{10}$, the number of singular points of $\bar H_1$ is equal
to $2^{16}$, and all its singular points are of the form
$w^4-z^4=0$.
\end{prop}
\proof A Hurwitz product $G(2,k)=G(2)^{\diamond k}$ is a Hurwitz
$C$-group of degree $m=2^{3k-1}$. By Propositions \ref{prod} and
\ref{2}, its Alexander polynomial is
$$\Delta _k(t)=(t-1)^k(t+1)^k.$$
By Theorem 6.2 from \cite{Ku}, each group $G(2,k)$ can be realized
as the fundamental group $\pi_1(\mathbb C^2\setminus H_k)$ for
some Hurwitz curve $\bar H_k$ of even degree and having
singularities of the form $w^{m}-z^{m}=0$. Since the multiplicity
of the root $t=1$ of the Alexander polynomial $\Delta _k(t)$ is
equal to $k$, it follows from Lemma \ref{number} and Theorem
\ref{mult} that the curve $\bar H_k$ consists of $k+1$ irreducible
components.

By Theorem \ref{th}, the first Betti number $b_1(\overline
X_{k,2})=k$, since the multiplicity of the root $t=-1$ of the
Alexander polynomial $\Delta _k(t)$ is equal to $k$.

To prove the existence of a Hurwitz curve $\bar H_1$ with the desired
properties, we should find some integer $m$ and a braid monodromy
factorization
$$\Delta_m^2=b_1\cdot \, .\, .\, .\, \cdot b_n$$
of a Hurwitz curve $\bar H$ such that the group given by the
presentation
$$<x_m,\dots,x_m\mid x_i=b_j(x_i)\,\,\text{for}\,\,
i=1,\dots,m,\,\, j=1,\dots,n>$$ is $C$-isomorphic to $G(2)$.

To find such presentation, let us recall briefly the proof of
Theorem 6.2 from \cite{Ku} and apply it to calculate the
invariants of a Hurwitz curve $\bar H_1$ for which $\pi_1(\mathbb
C\setminus H_1)\simeq G(2)$.

The group $G(2)$ is given by presentation (\ref{G2}). Consider the
group 
\begin{equation} \label{ya}
 G_{4,4} \simeq 
<x_1,\dots,x_4\mid x_i=\Delta_4^2(x_i)\,\, \text{for}\,\,
i=1,\dots,4>, 
\end{equation}
where $\Delta_4$ is the Garside element in $\text{Br}_4$. The
braid $\Delta_4^2$ is the braid monodromy of the singularity given by
$w^4-z^4=0$ and $s_0=\Delta^2_4$ (the factorization with a single
factor) is the braid monodromy factorization of four lines in
$\mathbb C\mathbb P^2$ passing through a fixed point.

To obtain presentation (\ref{G2}), we should add the relations
\begin{equation} \label{rele} x_2^2x_1x^{-2}_2=x_4,\, \, \,  x_4^2x_2x_4^{-2}=x_2,\, \,
x_3=x_2\end{equation} to presentation (\ref{ya}). For this, using
notations and notions from \cite{Ku}, in the beginning one should
perform the doubling (see Theorem 3.2 from \cite{Ku}) of the braid
monodromy factorization $s_0=\Delta_4^2$ several times in order to
have a possibility to move apart the generators $x_1,\dots,x_4$
and to change each relation from (\ref{rele}) by the relations
$x_i=x_{i+4}$ and relations of the form $x_i=b(x_i)$, where $b$ is
a braid conjugated to the standard generator of a braid group. It
is easy to see that in our case it suffices to perform the
doubling two times and we obtain the braid monodromy factorization
$s_1=d^2(s_0)$ of $\Delta^2_{16}$ (the doubling $d^2(s_0)$ is
defined in \cite{Ku} by formula (25)), each factor of which is
either conjugated to $\Delta_4^2$ or conjugated to a standard
generator of $\text{Br}_{16}$, and the group
\[
<x_1,\dots,x_{16}\mid x_i=b(x_i),\, \, i=1,\dots,16,
\, \, \text{and} \, \, b\,\, \text{is a factor of}\, \, s_1> 
\] is
$C$-isomorphic to $G_{4,4}$. The number of factors of $s_1$
conjugated to $\Delta_4^2$ is equal to $4^2$.

Then, to add the relations (\ref{rele}) to presentation (\ref{ya}),
one can use Lemma 3.4 from \cite{Ku} three times and obtain a
braid monodromy factorization $s_2$ of $\Delta^2_{2^{10}}$, each
factor of which is either conjugated to $\Delta_4^2$ or conjugated
to a standard generator of $\text{Br}_{2^{10}}$. The number of
factors of $s_2$ conjugated to $\Delta_4^2$ is equal to $4^8$. The
factorization $s_2$ is a braid monodromy factorization of a
Hurwitz curve $\bar H_1$, $\deg \bar H_1=2^{10}$, $\bar H_1$ has
$4^8$ singular points of the form $w^4-z^4=0$, and $\pi_1(\mathbb
C^2\setminus H)\cong G(2)$ by construction of $s_2$. \edvo

 \ifx\undefined\bysame
\newcommand{\bysame}{\leavevmode\hbox to3em{\hrulefill}\,}
\fi

\end{document}